\documentclass[11pt]{article}
\usepackage[english]{babel}
\usepackage[latin1]{inputenc}
\usepackage[T1]{fontenc}

\usepackage{amsmath, amsfonts, amssymb, amsthm}
\newtheorem{theo}{Theorem}
\newtheorem{lemma}{Lemma}
\newtheorem{prop}{Proposition}

\newtheorem{coro}{Corollary}

\newcommand{\R}{\mathbb{R}}
\newcommand{\Z}{\mathbb{Z}}
\newcommand{\N}{\mathbb{N}}
\newcommand{\E}{\mathbb{E}}
\renewcommand \P {\mathbb{P}}
\newcommand{\ind}{\mathbb{I}}

\newcommand{\bP}{\mathbf{P}}
\newcommand{\bE}{\mathbf{E}}

\newcommand{\la}{\langle}
\newcommand{\ra}{\rangle}
\newcommand{\rmP}{\mathrm{P}}
\newcommand{\rmE}{\mathrm{E}}
\newcommand{\ov}{\overline}
\newcommand \bOm {\mathbf{\Omega}}
\newcommand \G {\Gamma_{\varpi}}
\newcommand \vp {\varpi}
\newcommand \Ze {Z^{\vp,\e}}
\newcommand \Xe {X^{\vp,\e}}

\newcommand \Ee {\E_{[ \e^{-1}] \delta_0}}
\newcommand \Ge {\Gamma_{\vp}^{\e}}
\newcommand \kie {\chi^{\vp,\e}}
\newcommand \xie {\xi^{\vp,\e}}

\newcommand \Ue {U_{\vp,\e}}

\newcommand \cM {\mathcal{M}}
\newcommand \cK {\mathcal{K}}
\newcommand \cN {\mathcal{N}}
\newcommand{\cT}{\mathcal{T}}
\newcommand{\cR}{\mathcal{R}}
\newcommand{\cY}{\mathcal{Y}}

\newcommand{\e}{\varepsilon}
\newcommand \ka {\kappa}
\newcommand \vf {\varphi}

\def\build#1_#2^#3{\mathrel{
\mathop{\kern 0pt#1}\limits_{#2}^{#3}}}

\title{Escape probabilities for branching Brownian motion \\ among mild obstacles}
\author{Jean-Fran\c cois Le Gall, Amandine V\'eber}
\date{\small\today}

\begin{document}

\maketitle

\begin{abstract}
We derive asymptotics for the quenched probability that a critical branching Brownian motion
killed at a small rate $\e$ in Poissonian obstacles exits a large domain. Results
are formulated in terms of the solution to a semilinear partial differential equation with
singular boundary conditions. The proofs depend on a quenched homogenization 
theorem for branching Brownian motion among mild obstacles.

\smallskip
\noindent{\bf Keywords and phrases.} Branching Brownian motion, Poissonian obstacles, super-Brownian
motion, escape probability, homogenization, semilinear partial differential equation.

\smallskip
\noindent{\bf MSC2010 Classification Numbers.} Primary: 60K37, 60J80. Secondary: 60J68.
\end{abstract}

\section{Introduction}\label{intro BBM}

In the present work, we are interested in the long-term behaviour of branching Brownian motion killed in Poissonian obstacles. Let us start by describing a simple special case of our results. We consider a critical branching Brownian motion in $\R^d$ ($d\geq 1$), where all initial particles start from the origin. We assume that particles are killed at a (small) rate $\e>0$ within random balls of fixed radius, whose centers are distributed according to a homogeneous Poisson point process on $\R^d$. Then, how many initial particles do we need so that, with high probability, one of their descendants reaches distance $R$ from the origin ? Let $p_\e(R)$ be the (quenched) probability for our randomly killed branching Brownian motion starting with a single particle at $0$ to visit the complement of a large ball of radius $R$ centered at the origin. The preceding question is equivalent to determining the limiting behaviour of $p_\e(R)$ when $\e$ tends to $0$ and simultaneously $R$ tends to infinity.

The answer involves several regimes depending on the respective values of $\e$ and $R$. If $\e$ is small in comparison with $1/R^2$, the killing phenomenon does not matter and the result is the same as if there were no killing: $p_\e(R)$ behaves like a constant times $1/R^2$ (informally, the branching process must survive up to a time of order $R^2$ so that at least one of the particles travels a distance $R$, and well-known estimates for critical branching processes then lead to the correct asymptotics). On the other hand, if $\e$ is large in comparison with $1/R^2$, then the probability $p_\e(R)$ decreases exponentially fast as a function of $R\sqrt{\e}$: See Proposition \ref{large-dev} below.

Our main results focus on the critical regime where $\e R^2$ converges to a constant $a>0$. We show that the probability $p_\e(R)$ behaves like $R^{-2}$, as in the case without killing, but with a multiplicative constant which depends on $a$ and can be identified as the value at the origin of the solution of a semilinear partial differential equation with singular boundary conditions. A key tool to derive these asymptotics is a quenched homogenization theorem which shows that our branching Brownian motions among obstacles, suitably rescaled, are close to super-Brownian motion killed at a certain rate depending on $a$.

Let us formulate our assumptions more precisely in order to state our results. First, let us define the collection of obstacles. We denote the set of all compact subsets of $\R^d$ by $\cK$. This set is equipped with the usual Hausdorff metric $d_H$. Recall that $(\cK,d_H)$ is  a Polish space. For every $r>0$, $\cK_r$ denotes the subset of $\cK$ which consists of all compact sets that are contained in the closed ball of radius $r$ centered at the origin. Let $\Theta$ be a finite measure on $\cK$, and assume that $\Theta$ is supported on $\cK_{r_0}$ for some $r_0>0$. Let
$${\mathcal N}=\sum_{i\in I} \delta_{(x_i,K_i)}$$
be a Poisson point measure on $\R^d\times {\mathcal K}$ with intensity $\lambda_d\otimes \Theta$, where $\lambda_d$ stands for Lebesgue measure on $\R^d$. We assume that this point measure is defined on a probability space $({\bf\Omega},{\bf P})$ and we denote the generic element of ${\bf\Omega}$ by $\vp$. Our set of obstacles is then  defined by 
\begin{equation}
\label{obstacleset}
\Gamma_\varpi = \bigcup_{i\in I} (x_i + K_i),
\end{equation}
where obviously $x_i+ K_i = \{z=x_i + y : y\in K_i\}$. Note that we use the notation $\Gamma_\varpi $ to emphasize that the set of obstacles depends on the variable $\vp$ representing the environment. Let us also define a constant $\kappa$ by 
$$\kappa={\bf P}(0\in \Gamma_\vp)=1 - \exp\Big(-\int_{\mathcal K} \Theta(dK)\,\lambda_d(K)\Big).$$
To avoid trivial cases, we assume that $\kappa>0$, or equivalently $\Theta(\lambda_d(K)>0)>0$. By translation invariance, we also have ${\bf P}(x\in  \Gamma_\vp)=\kappa$ for every $x\in \R^d$.

Let us now introduce the sequence of branching Brownian motions of interest. Given $\vp\in {\bf \Omega}$ and a parameter $\e\geq 0$, we consider a branching Brownian motion on $\R^d$ such that
\begin{itemize}
\item each particle moves around in $\R^d$ according to the law of Brownian motion killed at rate $\e$ within  $\Gamma_\vp$ ;
\item each particle branches at rate $1$. During a branching event, the particle generates a random number of offspring, according to an offspring  distribution $\nu$ which has mean one and finite variance $\sigma^2\in(0,\infty)$.
\end{itemize}
This branching Brownian motion is denoted by $Z^{\vp,\e}=(Z^{\vp,\e}_t)_{t\geq 0}$, where $Z^{\vp,\e}_t$ stands for the sum of the Dirac point masses at the particles alive at time $t$. The processes $Z^{\vp,\e}$ are defined on a probability space $\Omega$. For every finite point measure $\mu$ on $\R^d$, we use the notation $\P_\mu$ for the probability measure on $\Omega$ under which each of the processes $Z^{\vp,\e}$ starts from $\mu$. 

Let $A$ be a bounded domain of class $C^2$ in $\R^d$ containing $0$.  We say that the branching Brownian motion $Z^{\vp,\e}$ hits $A^c$ if there exists $t>0$ such that $Z^{\vp,\e}_t(A^c)>0$. We are interested in asymptotics for the quantity 
$$\P_{\delta_0}(Z^{\vp,\e}\hbox{ hits } (RA)^c)$$
when $R\to\infty$ and $\e \to 0$. Here we use the obvious notation $RA=\{z=Ry:y\in A\}$.

\begin{theo}\label{theo charging proba}
For every $a\geq 0$, let $u_{(a)}=(u_{(a)}(x),x\in A)$ be the unique nonnegative solution of the singular boundary value problem
\begin{equation}
\label{BVpb}
\left\{\begin{array}{ll}
{\displaystyle\frac{1}{2}\Delta u = \frac{\sigma^2}{2}\,u^2 +  a\,u}\qquad &\hbox{in }A,\\
\noalign{\smallskip}
u_{|\partial A} = +\infty\;.&
\end{array}
\right.
\end{equation}
Then,
$$\lim_{R\to\infty}\Big(\sup_{\e\geq 0}\Big|R^2 \P_{\delta_0}(Z^{\vp,\e}\hbox{ hits } (RA)^c) -u_{(\kappa \e R^2)}(0)\Big|\Big) = 0\;,\qquad {\bf P}(d\vp)\ \hbox{a.s.}$$
\end{theo}

Let us state a corollary of the theorem, which is motivated by the simple question we asked
at the beginning of this introduction. For every $a\geq 0$, we denote by $u^\circ_{(a)}$ the solution of the 
boundary problem (\ref{BVpb}) when $A$ is the open unit ball of $\R^d$. 

\begin{coro}
\label{escapepro}For every $\e\in(0,1)$, let $n_\e$ be a positive integer. Assume that $\e n_\e\longrightarrow b$ as $\e\to 0$, for some $b>0$. Denote by $R^{\varpi,\e}$ the maximal distance from the origin
attained by a particle of the branching Brownian motion $Z^{\vp,\e}$. Then,
$\bP(d\varpi)$ a.s., the law of $\sqrt{\e}\,R^{\varpi,\e}$ under $\P_{n_\e\delta_0}$
converges as $\e\to 0$ towards the probability measure $\pi_b$ on $\R_+$ defined by
$$\pi_b([0,r])=\exp(-\frac{b}{r^2}\,u^\circ_{(\kappa  r)}(0))$$
for every $r>0$.
\end{coro}

In the setting of Theorem \ref{theo charging proba},
it is not hard to see that $u_{(a)}(0)$ tends to $0$ as $a\to\infty$ (see Lemma \ref{analytic} below)
and thus this theorem does not give much information when $R\to\infty$ and $\e\to 0$ in such a way that $\e R^2$ tends to $\infty$. In that case, the next proposition provides an exponential decay, which contrasts with the preceding theorem. Since our bounds are clearly not optimal, we consider only
the case when $A$ is a ball. We denote the open ball of radius $r$ centered at the origin
by $B(0,r)$. In the general case we may apply the bounds (i) and (ii) of the proposition after replacing $A$ by a ball $B(0,r)$ such that $B(0,r)\supset A$ or $B(0,r)\subset A$ respectively. 

\begin{prop}
\label{large-dev}
{\rm (i)} There exists a positive constant $C_0=C_0(\nu)$ such that, for every $\varpi\in {\bf\Omega}$, $R\geq 1$ and $\e\in[1/R^2,1]$,
$$\P_{\delta_0}(Z^{\vp,\e} \hbox{ hits } B(0,R)^c)\geq C_0\e\,\exp(-R\sqrt{2\e}).$$
{\rm (ii)} There exists a positive constant $C_1=C_1(\Theta)$ such that for every $R\geq 1$ and $\e\in[1/R^2,1]$,
$${\bf P} \otimes \P_{\delta_0}(Z^{\vp,\e} \hbox{ hits } B(0,R)^c)\leq \exp(-C_1 R\sqrt{\e}).$$
Consequently, we can find two positive constants $C_2=C_2(\Theta)$ 
and $C_3=C_3(\Theta)$ such that $\bP(d\varpi)$ a.s., for every
sufficiently large $R$ and every $\e\in [C_2(\log\log R)^2/R^2, 1]$,
$$\P_{\delta_0}(Z^{\vp,\e} \hbox{ hits } B(0,R)^c)\leq \exp(-C_3 R\sqrt{\e}).$$
\end{prop}

Part (i) of the proposition is derived from an estimate about branching Brownian motion killed homogeneously at rate $\e$,
which explains why this bound holds for every $\varpi\in {\bf\Omega}$ and does not depend on the measure $\Theta$. The bounds in (ii) follow from an estimate for Brownian motion 
killed in mild obstacles and therefore do not depend on the offspring distribution $\nu$. The first assertion in (ii) may be compared to Proposition 5.2.8 
in \cite{SZN1998}.

Proposition \ref{large-dev} only gives rather crude estimates, and it would be of interest to obtain more precise information on the decay of the quenched probabilities $\P_{\delta_0}(Z^{\vp,\e}\hbox{ hits } (RA)^c)$ in the case when $\e R^2$ tends to $\infty$. This leads to large deviations problems in the spirit of the work of Sznitman \cite{SZN1998}, which we do not address here.

A major ingredient of the proof of Theorem \ref{theo charging proba} is the following quenched homogenization result. We need to introduce a rescaled version of the process $Z^{\vp,\e}$. For every $\e>0$ and every $t\geq 0$, let us define a random measure $X^{\vp,\e}_t$ on $\R^d$ by setting,
for any nonnegative measurable function $\varphi$ on $\R^d$,
\begin{equation}
\label{scalingtransfo}
\langle X^{\vp,\e}_t,\varphi\rangle = \e \int Z^{\vp,\e}_{\e^{-1}t}(dx)\,\varphi(\e^{1/2} x).
\end{equation}
Here and later, the notation $\langle \mu, \varphi\rangle$ stands for the integral of the
function $\varphi$ against the measure $\mu$, whenever this integral makes sense.

For every real $x\geq 0$, $[x]$ denotes the integer part of $x$.

\begin{theo}\label{theo cv to sbm}
Except for a $\bP$-negligible set of values of $\vp$, the law of $(X^{\vp,\e}_t)_{t\ge0}$ under $\P_{[\e^{-1}]\delta_0}$ converges weakly as $\e\to 0$, in the Skorokhod sense, to that of a super-Brownian motion with branching mechanism $\psi_{(\kappa)}(u):=\frac{\sigma^2}{2}u^2 + \kappa u$ started at $\delta_0$. 
\end{theo} 
The definition of super-Brownian motion with branching mechanism $\psi_{(\kappa)}$ is
recalled in Section 2 below. 

As a hint of why Theorem \ref{theo cv to sbm} should be true, notice that for a given realization of the obstacles, the probability that a single Brownian motion
starting from $0$ and killed at rate $\e$ within $\G$ is still alive by time $t>0$ is given by
\begin{equation}\label{intro conv}
\rmE\Big[\exp\Big\{-\e \int_0^t \ind_{\G}(\xi_s)ds\Big\}\Big],
\end{equation}
where $\xi$ denotes standard $d$-dimensional Brownian motion. Let us focus on the integral within the exponential in (\ref{intro conv}). Averaging over the law of the obstacles and using Fubini's theorem, we obtain for each $t\geq 0$
$$
\bE\Big[\rmE\Big[\e \int_0^t \ind_{\G}(\xi_s)ds\Big]\Big]= \e \ \rmE\Big[\int_0^t\bP[\xi_s\in \G]ds\Big]=\e \ka t.
$$
We can thus guess,
and easily prove, that the rescaled Brownian motion $(\sqrt{\e}\ \xi_{\e^{-1}t},\ t\geq 0)$, 
which is killed at rate $1$ within $\sqrt{\e}\ \G$, converges to Brownian motion killed at homogeneous rate $\ka$ as $\e\rightarrow 0$. Theorem \ref{theo cv to sbm} shows that an analogous convergence indeed holds in our more general framework of branching Brownian motions, for any fixed $\vp$ contained in a set of $\bP$-probability one.

Let us briefly explain how Theorem \ref{theo charging proba} is derived from Theorem \ref{theo cv to sbm}.
Consider a sequence $(\e_n,R_n)$ such that $R_n\to \infty$ and $\e_nR_n^2$ 
converges to a positive constant $a$. By a simple scaling transformation, the probability
$\P_{\delta_0}(Z^{\vp,\e_n}\hbox{ hits } (R_nA)^c)$ coincides with 
$\P_{\delta_0}(X^{\vp,\e_n}\hbox{ hits } (b_nA)^c)$, where $b_n=\e_n^{1/2}R_n$ converges
to $\sqrt{a}$. We can then use Theorem \ref{theo cv to sbm} to investigate the asymptotic behaviour
of the latter
hitting probabilities. This limiting behaviour involves the corresponding hitting probabilities 
for super-Brownian motion, which are known to be related to solutions of
semilinear partial differential equations from the work of Dynkin \cite{DYN1991,DYN1993}. One difficulty in
implementing the preceding idea comes from the fact that the convergence in Theorem \ref{theo cv to sbm} is not strong enough to ensure that hitting probabilities for the processes 
$(X^{\vp,\e}_t)_{t\ge0}$ converge to hitting probabilities for the limiting process. Much of the proof
of Theorem \ref{theo charging proba} in Section 4 is devoted to a precise justification of this property (Lemma
\ref{convhittingpro}).

To complete this introduction, let us mention that branching Brownian motion and superprocesses among random obstacles
have been studied recently in several papers, including Engl\"ander and den Hollander \cite{EDH2003} and 
Engl\"ander \cite{ENG2008}. These papers concentrate on the case of supercritical branching, in contrast with
critical branching which is considered here. See also the survey \cite{ENG2007}. A homogenization theorem 
related to Theorem \ref{theo cv to sbm} has been proved in \cite{VEB2009} for super-Brownian motion among hard obstacles,
in the case when the intensity of the obstacles grows to infinity but their diameters shrink to $0$. There is a huge
literature about Brownian motion and random walks among (hard or mild) obstacles, and the reader may look at the book of Sznitman \cite{SZN1998} for additional references.

The rest of this paper is laid out as follows. In Section \ref{section: preliminaries}, we introduce the basic notation and objects, and state several results about hitting probabilities for spatial branching processes we shall need in the sequel. Theorem \ref{theo cv to sbm} 
and Proposition \ref{large-dev} are proved in Section \ref{section: quenched conv}. Theorem \ref{theo charging proba} and Corollary \ref{escapepro} are then derived in Section \ref{section: main result}.
\section{Preliminaries}\label{section: preliminaries}

\subsection{Notation}

We denote the set of all finite measures on $\R^d$ by $\cM_f(\R^d)$. This set is equipped with the weak topology. We write $\cM_p(\R^d)$ for the subset 
of ${\cal M}_f(\R^d)$ which consists of all finite point measures on $\R^d$. 

If $E$ is a metric space, we denote  the set of all bounded continuous functions on the space $E$ by $\bar{C}(E)$ and we let $\|f\|$ stand for the supremum norm of $f\in \bar{C}(E)$. We write $C^2(\R^d)$ for the set of all twice continuously differentiable functions on $\R^d$, and $\bar{C}^2(\R^d)$ for that of all bounded functions in $C^2(\R^d)$ whose first and second derivatives are also bounded. An index $+$ added to this notation means that we require the functions to be nonnegative.  We equip $\bar{C}^2(\R^d)$ with the topology induced by the seminorms $\|f\|_{(R)}$, where for every $R>0$ 
$$
\|f\|_{(R)}:= \sup_{|x|\leq R}\bigg\{|f(x)|+\sum_{i=1}^d \Big|\frac{\partial f}{\partial x_i}(x)\Big|+ \sum_{i,j=1}^d \Big|\frac{\partial^2 f}{\partial x_i \partial x_j}(x)\Big| \bigg\}.
$$

If $E$ is a Polish space, we let $D_{E}[0,\infty)$ be the set of all c\`adl\`ag paths with values in $E$, equipped with the Skorokhod topology.

If $x\in\R^d$ and $r>0$, $B(x,r)$ denotes the open ball of radius $r$ centered at $x$, and $\bar B(x,r)$ stands for the
corresponding closed ball. More generally, the closure of a subset $F$ of $\R^d$ is denoted by $\bar F$. Lebesgue measure on $\R^d$
is denoted by $\lambda_d$. 

Finally, the notation $\xi=(\xi_t)_{t\geq 0}$ will stand for a standard Brownian motion in $\R^d$, which starts from $x$
under the probability measure $\rmP_x$. It will also be convenient to use the notation 
$\xi^{\varpi,\e}$ for Brownian motion in $\R^d$ killed at rate $\e$ in the set $\Gamma_\varpi$.
As usual, the value of $\xi^{\varpi,\e}$ after its killing time is a cemetery point $\Delta$
added to $\R^d$, and we agree that all functions vanish at $\Delta$. 

\subsection{Super-Brownian motion}

Let $a\geq 0$ and set $\psi_{(a)}(u)=\frac{\sigma^2}{2} u^2 + a\,u$, for every $u\geq 0$ (the offspring
distribution $\nu$, and thus the parameter $\sigma>0$ are fixed throughout this work). Super-Brownian motion with branching mechanism $\psi_{(a)}$ is the continuous strong Markov process with values in $\cM_f(\R^d)$, whose transition kernels $(Q_t)_{t\geq 0}$ are characterized as follows: For every $g\in \bar{C}_+(\R^d)$ and every $\mu\in \cM_f(\R^d)$, we have for every $t\geq 0$
\begin{equation}
\label{transiker}
\int Q_t(\mu,d\mu')\,\exp(-\la \mu', g\ra) = \exp(-\la \mu, V_tg\ra),
\end{equation}
where the function $u_t(x)=V_tg(x)$, $t\geq 0,\ x\in\R^d$, is the unique nonnegative solution of the semilinear parabolic problem
$$\left\{ \begin{array}{l}
{\displaystyle\frac{\partial u}{\partial t} }= \frac{1}{2} \Delta u - \psi_{(a)}(u)\quad\hbox{ in }(0,\infty)\times \R^d\;,\\
\noalign{\smallskip}
u_0=g\;.
\end{array}
\right.
$$

Let $Y=(Y_t)_{t\geq 0}$ be a super-Brownian motion with branching mechanism $\psi_{(a)}$, started at $\mu\in \cM_f(\R^d)$. Then, for every $g\in \bar C_+^2(\R^d)$
\begin{equation}\label{limit mp}
e^{-\la Y_t,g\ra} -e^{-\la Y_0,g\ra}-\int_0^t \Big\la Y_s,-\frac{1}{2}\Delta g +\psi_{(a)}(g)\Big\ra\ e^{-\la Y_s, g\ra}ds
\end{equation}
is a martingale. It is well known that this martingale problem and the initial value $\mu$ characterize the law of $Y$. This is indeed an application of the classical ``duality method'' (see in particular Chapter 4 in \cite{EK1986}). The nonlinear semigroup $g\rightarrow V_tg$ provides a deterministic dual to super-Brownian motion, and the duality argument then shows that if a measure-valued process started from $\mu$ satisfies the preceding martingale problem, the Laplace functional of its value at time $t$ must be given by the right-hand side of (\ref{transiker}). See Section 1.6 of \cite{ETH2000} for more details.

\subsection{Branching Brownian motion among random obstacles}\label{subs: BBM et ob}

In view of our applications (and in particular because we want to refer to some results of \cite{Ch}), it will be convenient to give a more formal description of the branching Brownian motions that were already introduced in Section \ref{intro BBM} above. Recall that our offspring distribution $\nu$ is assumed to be critical and that $\mathrm{var}(\nu)=\sigma^2\in(0,\infty)$. The probability generating function of  $\nu$ 
will be denoted by $\Upsilon$. 

Let $\cT$ be a Galton-Watson tree with offspring distribution $\nu$ (see e.g. \cite{LEG2005}). As usual, we view $\cT$ as a random finite subset of $$U:=\bigcup_{n=0}^{\infty} \N^n,$$ 
where $\N=\{1,2,3,\ldots\}$, and 
$\N^0=\{\varnothing\}$. If $v=(v_1,\ldots,v_{n})\in U\backslash \{\varnothing\}$, the parent of $v$ is denoted by $\hat v=(v_1,\ldots,v_{n-1})$ and we also use the notation $v\prec v'$ to mean that $v'$ is a descendant of $v$ distinct from $v$. Consider a collection $(e_v,v\in U)$ of independent exponential random variables with parameter $1$, which is also independent of $\cT$. We define for every $v\in U$ its birth time $\alpha_v$ and its death time $\beta_v$ recursively by setting $\alpha_\varnothing=0$ and  $\beta_\varnothing = e_\varnothing$, and for every $v\in U\backslash\{\varnothing\}$,
$$\alpha_v = \beta_{\hat v}\ ,\ \beta_v=\alpha_v + e_v.$$
Let us now construct the spatial motions. Fix a starting point $x\in \R^d$, and consider a collection $(B^v,v\in U)$ of independent standard Brownian motions in $\R^d$ (started from $0$), independent of $\cT$ and of $(e_v,v\in U)$. For every $v\in U$, define the \emph{historical path} $\omega^v=(\omega^v_t,0\leq t\leq \beta_v)$ associated with $v$ in the following way. First  $\omega^\varnothing_t=x+B^\varnothing_t$ for $0\leq t\leq \beta_\varnothing$. Then, if $v\in U\backslash\{\varnothing\}$, set $\omega^v_t=\omega^{\hat v}_t$ for all $0\leq t\leq\alpha_v$ and
$$\omega^v_t= \omega^{\hat v}_{\alpha_v} + B^v_{t-\alpha_v}\qquad\hbox{for }
\alpha_v\leq t\leq \beta_v.$$
A branching Brownian motion (without killing in obstacles) starting from $\delta_x$ is obtained by setting for every $t\geq 0$, 
$$Z_t=\sum_{v\in\cT,v\sim t} \delta_{\omega^v_t},$$
where the notation $v\sim t$ means that $\alpha_v\leq t< \beta_v$.

In this formalism, it is now easy to introduce killing in obstacles. Consider yet another independent collection $(\gamma_v)_{v\in U}$ of independent exponential random variables with parameter $1$, and define for every $\vp\in\bOm$, $\e\geq 0$, and for every $v\in U$
$$\zeta^{\vp,\e}_v:= \inf\bigg\{s\in [\alpha_v,\beta_v): \int_{\alpha_v}^s dr\,\ind_{\Gamma_\vp} (\omega^v_r) > \e^{-1}\,\gamma_v\bigg\},$$
where $\inf\varnothing =\infty$. By setting
$$Z^{\vp,\e}_t=\sum_{v\in\cT,v\sim t}  \ind_{\{t<\zeta^{\vp,\e}_v\;{\rm and}\; \zeta^{\vp,\e}_{v'}=\infty, \;{\rm for}\; {\rm every}\; v'\prec v\}}\, \delta_{\omega^v_t},$$
we obtain a branching Brownian motion killed at rate $\e$ in the obstacle set $\Gamma_{\vp}$, starting from $\delta_x$. An obvious extension of the preceding construction allows us to obtain branching Brownian motions starting from any point measure $\mu\in \cM_p(\R^d)$. 

We now recall a special case of the classical convergence of rescaled branching Brownian motions towards super-Brownian motion. For our applications, we state the case where particles are killed at a constant rate homogeneously over $\R^d$ (this case is obtained from the preceding construction of $Z^{\vp,\e}$ by replacing $\G$ by $\R^d$). In the next two statements, for every $\e>0$, $Z^{(\e)}=(Z^{(\e)}_t)_{t\geq 0}$ denotes a branching Brownian motion with offspring distribution $\nu$, where particles are killed homogeneously over $\R^d$ at rate $\e$. As previously $Z^{(\e)}$
starts from $\mu$ under the probability measure $\P_\mu$, for every $\mu\in \cM_p(\R^d)$.

\begin{prop}
\label{convBBM}
Let $a\geq 0$. For every $\e>0$, define a measure-valued process $(X^{(\e)}_t)_{t\geq 0}$ by setting 
$$\langle X^{(\e)}_t,\vf\rangle  =\e\int Z^{(a\e)}_{\e^{-1}t}(dx)\,\vf(\e^{1/2}x),$$
for every $\vf\in\bar C(\R^d)$. For every fixed $\eta>0$,  the law of $X^{(\e)}$ under  
$\P_{[\eta\e^{-1}]\delta_0}$ converges as $\e\to 0$, in the Skorokhod sense,  
towards the law of super-Brownian motion with branching mechanism $\psi_{(a)}$ starting from $\eta\delta_0$.
\end{prop}
A proof of Proposition \ref{convBBM} can be found in Chapter 1 of \cite{ETH2000} in the case $a=0$,
and arguments are easily adapted to cover the general case.

Finally, we shall use an estimate  for the probability that a branching Brownian motion starting from $\delta_0$ exits a large ball centered at the origin. 
Similar estimates can be found in Sawyer and Fleischman \cite{SF1979}, but we provide a short proof for the sake of completeness.
\begin{lemma}
\label{Sawyer}
Suppose that $d=1$.
There exist two positive constants $C'_0=C'_0(\nu)$ and $C'_1=C'_1(\nu)$ such that, for every $\e\in[0,1]$ and  $r\geq 1$,
$$C'_0( r^{-2}\,\ind_{\{r\leq \frac{1}{\sqrt{\e}}\}}+ \e e^{-r\sqrt{2\e}}
\,\ind_{\{r> \frac{1}{\sqrt{\e}}\}})\leq \P_{\delta_0}( Z^{(\e)} \hbox{ hits } (-r,r)^c) \leq C'_1( r^{-2}\,\ind_{\{r\leq \frac{1}{\sqrt{\e}}\}}+ \e e^{-r\sqrt{2\e}}
\,\ind_{\{r> \frac{1}{\sqrt{\e}}\}}) .$$
\end{lemma}

\noindent{\bf Remark}. As an immediate consequence of the upper bound of the lemma, we have 
in dimension $d$, for every $r>0$,
\begin{equation}
\label{cons-Sawyer}
\P_{\delta_0}( Z^{(0)} \hbox{ hits } B(0,r)^c) \leq C''_1\,(r+1)^{-2}
\end{equation}
with a constant $C''_1=C''_1(d,\nu)$.

\medskip

\noindent{\bf Proof}. It clearly suffices to prove that the stated bounds hold for
the quantity $\P_{\delta_r}( Z^{(\e)} \hbox{ hits } (-\infty,0])$
instead of  $\P_{\delta_0}( Z^{(\e)} \hbox{ hits } (-r,r)^c)$. We fix $\e\geq 0$, and for every 
$x> 0$ and $t\geq 0$, we set
$$q_\e(x,t)= \P_{\delta_x}\left(Z^{(\e)}\hbox{ does not hit } (-\infty,0]\hbox{ before time }t\right).$$
and 
$$p_\e(x)= \P_{\delta_x}\left(Z^{(\e)}\hbox{ hits }(-\infty,0]\right)= \lim_{t\uparrow\infty} \uparrow (1- q(x,t)).$$
In this proof only, we write $\rmP^\e_x$
for the probability under which $\xi$ is a Brownian motion starting from $x$ and killed at rate $\e$
(upon killing, $\xi$ is sent to the cemetery point $\Delta$ and we recall that all
functions vanish at $\Delta$).
Write $S:=\inf\{t\geq 0: \xi_t\in(-\infty,0]\}$. By standard arguments (see e.g. the proof of Proposition
II.3 in \cite{LEG1999}), the function $q_\e$ solves the integral equation
$$q_\e(x,t)=\rmP^\e_x(S>t)+ \rmE^\e_x\Big[\int_0^{t\wedge S} (\Upsilon(q_\e(\xi_s,t-s))-q_\e(\xi_s,t-s))\,ds\Big],$$
where we recall that $\Upsilon$ denotes the generating function of the offspring distribution $\nu$.
For every $a\in [0,1]$, set $\Phi(a)=\Upsilon(1-a)- (1-a)$. Note that $\Phi(0)=0$ and the function $\Phi$
is monotone increasing under our assumptions.
Furthermore, $\Phi(a)=\frac{\sigma^2}{2}a^2 +o(a^2)$ when $a\to 0$. By a monotone passage to the limit we get that, for every
$x>0$,
\begin{equation}
\label{intequ}
p_\e(x)+ \rmE^\e_x\Big[\int_0^{S} \Phi(p_\e(\xi_s))\,ds\Big] = \rmP_x^\e(S<\infty).
\end{equation}
It follows that the function $p_\e$ satisfies the differential equation
$$\frac{1}{2}p_\e'' = \e\,p_\e + \Phi(p_\e)$$
on $(0,\infty)$ with boundary conditions $p_\e(0)=1,p_\e(\infty)=0$. By solving this differential equation, we get, for every $x>0$,
$$\int_{p_\e(x)}^1 \frac{du}{\sqrt{2\e u^2 + 4\Gamma(u)}} = x$$
where $\Gamma(u)=\int_0^u \Phi(v)\,dv$. Note that there exist positive constants $c,c'$ such that
$cu^3\leq\Gamma(u)\leq c'u^3$ for every $u\in[0,1]$. The desired bounds then follow from easy analytic arguments.
$\hfill\square$

\subsection{Hitting probabilities for super-Brownian motion}\label{subs: hitting prob}

Let $Y^{(a)}=(Y^{(a)}_t)_{t\geq 0}$ be a super-Brownian motion with branching mechanism $\psi_{(a)}$ for some $a\geq 0$. Suppose that $Y^{(a)}$ starts from $\mu$ under the probability measure $P_{\mu}$, for every $\mu\in\cM_f(\R^d)$.

The range of $Y^{(a)}$ is by definition
$$\cR(Y^{(a)})= \bigcup_{\e>0}\left(\overline{\bigcup_{t=\e}^\infty {\rm supp}(Y^{(a)}_t)}\right),$$
where for every $\mu\in\cM_f(\R^d)$, ${\rm supp}(\mu)$ denotes the topological support of $\mu$. 

Let $D$ be a domain in $\R^d$ and let $x\in D$. Consider the process $Y^{(a)}$ started from $\delta_x$. We say that $Y^{(a)}$ hits $D^c$ if the range $\cR(Y^{(a)})$ intersects $D^c$. By a famous result of Dynkin \cite{DYN1991,DYN1993} the function
$$u^D_{(a)}(x) = -\log\,\Big(1-P_{\delta_x}(\cR(Y^{(a)})\cap D^c\not =\emptyset)\Big)\;,\quad x\in D\,,$$
is the maximal nonnegative solution of the semilinear partial differential equation $\frac{1}{2}\Delta u = \psi_{(a)}(u)$ in $D$. 

Under mild regularity assumptions on $D$ (which hold e.g. when $D$ satisfies an exterior cone condition at every point of $\partial D$), the function $u^D_{(a)}$ has boundary value $+\infty$ at every point of $\partial D$ and is the unique nonnegative solution of the equation $\frac{1}{2}\Delta u = \psi_{(a)}(u)$ in $D$ with boundary value $+\infty$ everywhere on $\partial D$. A discussion of this result and related ones can be found in 
Chapter VI of the book \cite{LEG1999}. This reference considers only the case $a=0$, but the same results can be obtained for any 
$a\geq 0$ by similar arguments: Note that the Brownian snake approach can be extended from the case $a=0$ considered
in \cite{LEG1999} to $a\geq 0$, simply by replacing the reflecting Brownian motion driving the snake by a reflecting Brownian
motion with negative drift (see Chapter 4 of \cite{DLG2002} for a discussion of the snake approach to
superprocesses with a general branching mechanism).

We shall be interested in the special case $D=A$. Recall that $0\in A$ and that we assume $A$ is a domain of class $C^2$, meaning that the boundary of $A$ can be represented locally as the graph of a twice continuously differentiable function, in a suitable  system of coordinates. We write $u_{(a)}(x)=u^A_{(a)}(x)$ to simplify notation. From the analytic viewpoint, the function $u_{(a)}$ may be constructed as follows. For every integer $n\geq 1$, let $u_{(a),n}$ be the unique nonnegative solution of the nonlinear Dirichlet problem 
$$\left\{\begin{array}{l}
\frac{1}{2}\Delta u = \psi_{(a)}(u)\quad\hbox{in }A\,,\\
u_{|\partial A} = n\,.
\end{array}
\right.
$$
Then $u_{(a)}=\lim\uparrow u_{(a),n}$ as $n\to\infty$.

The following lemma records certain analytic properties which will be useful in the forthcoming proofs.

\begin{lemma}
\label{analytic}
{\rm (i)} Let $x\in A$. The function $a\longrightarrow u_{(a)}(x)$ is continuous 
and nonincreasing on $[0,\infty)$, and tends to $0$ as $a\to \infty$.
\par\noindent{\rm (ii)} For every $\delta\in (0,{\rm dist}(0,A^c))$, let $A_\delta$ be the subdomain of $A$ defined as the connected component of the open set $\{x\in A: {\rm dist}(x,A^c)>\delta\}$ that contains $0$. Then, for every $a\geq 0$, $u^{A_\delta}_{(a)}(0)$ tends to $u_{(a)}(0)$ as $\delta \to 0$. 
\end{lemma}

\noindent{\bf Proof}. (i) Let us first verify that the function $a\longrightarrow u_{(a)}(x)$ is monotone nonincreasing, for every $x\in A$. To see this, we apply a standard comparison principle (see e.g. Lemma V.7 in \cite{LEG1999}) to obtain that $u_{(a'),n}\leq u_{(a),n}$ if $a\leq a'$, for every $n\geq 1$. It then suffices to let $n\to\infty$. 

Let $(a_k)_{k\geq 1}$ be a sequence of nonnegative reals increasing to $a\in(0,\infty)$. We can set for every $x\in A$ 
$$v(x)=\lim_{k\uparrow \infty}\downarrow u_{(a_k)}(x),$$
and we have $v\geq u_{(a)}$. In order to verify that $v\leq u_{(a)}$, we only need to check that $v$ solves $\frac{1}{2}\Delta v = \psi_{(a)}(v)$ in $A$ (recall that $u_{(a)}$ is the maximal nonnegative solution of this equation). To do so, let $B$ be an open ball whose closure $\bar B$ is contained in $A$. For every $k\geq 1$, the restriction of $u_{(a_k)}$ to $B$ solves the equation $\frac{1}{2}\Delta u = \psi_{(a_k)}(u)$ in $B$. By the probabilistic interpretation of the integral equation associated with this PDE (see e.g. Chapter V in \cite{LEG1999}), this implies that, for every $x\in B$, 
$$u_{(a_k)}(x)+\rmE_x\bigg[\int_0^{\tau_B} \psi_{(a_k)}(u_{(a_k)}(\xi_s))\,ds \bigg] = \rmE_x\big[u_{(a_k)}(\xi_{\tau_B})\big],$$
where we recall our notation $\xi$ for a Brownian motion starting from $x$ under the probability measure $\rmP_x$, and $\tau_B:=\inf\{t\geq 0:\xi_t\notin B\}$. By passing to the limit $k\to\infty$ in the previous display, we can write
$$v(x)+\rmE_x\bigg[\int_0^{\tau_B} \psi_{(a)}(v(\xi_s))\,ds \bigg] = \rmE_x\big[v(\xi_{\tau_B})\big],$$
which is enough to obtain that $v$ solves $\frac{1}{2}\Delta v = \psi_{(a)}(v)$ in $B$, and therefore in $A$ since $B$ was arbitrary.

Similar arguments show that, if $(a_k)_{k\geq 1}$ is a decreasing sequence of nonnegative reals converging to $a\in[0,\infty)$, then $u_{(a_k)}(x)$ converges to $u_{(a)}(x)$ for every $x\in A$. Finally, the fact that $u_{(a)}(x)$ tends to $0$ as $a\to\infty$ can be obtained from the comparison principle: If $B$ is a ball such that $\bar B\subset A$, the restriction of $u_{(a)}$ to $B$ is bounded above by the solution $v_{(a)}$ of the linear equation $\frac{1}{2}\Delta v_{(a)}=a\,v_{(a)}$ in $B$, with boundary value equal to the restriction of $u_{(0)}$ to $B$. It is easily seen that $v_{(a)}(x)\longrightarrow 0$ as $a\to\infty$, for instance by using the Feynman-Kac formula.

\noindent (ii) Fix $a\geq 0$. If $0<\delta<\delta'$, the closure of $A_{\delta'}$ is contained in $A_{\delta}$. The restriction of $u^{A_{\delta}}_{(a)}$ to $A_{\delta'}$ is a nonnegative solution of  $\frac{1}{2}\Delta u = \psi_{(a)}(u)$ in $A_{\delta'}$ and is thus bounded above by $u_{(a)}^{A_{\delta'}}$. Hence, for every fixed $x\in A$ the function $\delta \to u^{A_\delta}_{(a)}(x)$, which is defined for $\delta>0$ small enough, is nondecreasing and we can set
$$v(x)=\lim_{\delta\downarrow 0} u^{A_\delta}_{(a)}(x).$$
By the same argument we used to obtain the monotonicity of the mapping $\delta \to u^{A_\delta}_{(a)}(x)$, we also have $v(x)\geq u_{(a)}(x)$ for every $x\in A$. To obtain the reverse inequality $v\leq u_{(a)}$, it is enough to verify that $v$ solves $\frac{1}{2}\Delta v = \psi_{(a)}(v)$ in $B$. But this follows by arguments similar to those we used in the proof of part (i) of the lemma. $\hfill\square$

\begin{lemma}
\label{hittinglemma}
For every $a\geq 0$ and $x\in A$,
$$\{\cR(Y^{(a)})\cap A^c \not =\emptyset\} = \{\cR(Y^{(a)})\cap \bar A^c\not =\emptyset\}\,,\quad P_{\delta_x}\ {\rm a.s.}$$
\end{lemma}

\noindent{\bf Proof}. The inclusion
$$\{\cR(Y^{(a)})\cap A^c \not =\emptyset\} \supset \{\cR(Y^{(a)})\cap \bar A^c\not =\emptyset\}$$
is trivial. To show the reverse inclusion, we may argue as follows. By Theorem IV.9 in \cite{LEG1999} (which holds under much less stringent assumptions on $A$), the event  $\{\cR(Y^{(a)})\cap A^c \not =\emptyset\}$ holds if and only if the exit measure of the super-Brownian motion $Y^{(a)}$ from $A$ is nonzero. Applying the special Markov property of superprocesses
\cite[Theorem 1.3]{DYN1993}, we see that it is enough to prove that for super-Brownian motion starting from a nonzero initial measure supported  on $\partial A$, the range immediately hits $(\bar A)^c$. This is however easy under our regularity assumptions on $A$. We leave the details to the reader. $\hfill\square$

\smallskip 
From now on, we write $\{Y^{(a)}\hbox{ hits } F\}$ for the event $\{\cR(Y^{(a)})\cap F\not =\emptyset\}$.

\section{Quenched convergence to super-Brownian motion}\label{section: quenched conv}

The main goal of this section is to prove Theorem \ref{theo cv to sbm}. At the end of the section, we also
establish Proposition \ref{large-dev}, using certain arguments related to the proof of
Theorem  \ref{theo cv to sbm}. To simplify notation, we set for every $\vp\in \bOm$ and $\e\in (0,1)$,
$$
\Ge = \sqrt{\e} \ \G = \bigcup_{i\in I}\ \sqrt{\e}(x_i+K_i).
$$

The following lemma identifies a martingale problem solved by our branching Brownian motion $\Ze$. It can be proved by standard arguments (see e.g. Section 9.4 in \cite{EK1986}).

\begin{lemma}\label{lemma: mart pb} Let $\vp\in \bOm$ and $\e\geq 0$. Under each probability $\P_\mu$, $\mu\in {\cal M}_p(\R^d)$, the process $\Ze$ solves the following martingale problem: For every $f\in \bar C_+^2(\R^d)$ such that $0<\inf f\leq f\leq 1$, the process
$$
e^{\la \Ze_t,\log f\ra}-e^{\la \Ze_0,\log f\ra}-\int_0^t \bigg\la \Ze_s,\frac{\frac{1}{2}\Delta f+\e \ind_{\G}(1-f)+\Upsilon(f)-f}{f} \bigg\ra e^{\la \Ze_s,\log f\ra} ds
$$
is a martingale.
\end{lemma}

We can derive from Lemma \ref{lemma: mart pb} (or from a direct argument) that for every $g\in \bar C^2(\R^d)$, 
\begin{equation}\label{mp Xe-bv}
M_t(g):= \la \Ze_t,g\ra - \la \Ze_0,g \ra - \int_0^t \Big\la \Ze_s,\frac{1}{2}\Delta g -\e \ind_{\G}g\Big\ra\ ds
\end{equation}
is a martingale. An easy computation gives that the square bracket of this martingale is
\begin{equation}\label{mp Xe-qv}
[M(g),M(g)]_t= \int_0^t \la \Ze_s,\nabla g.\nabla g \ra \ ds + \sum_{0\leq s\leq t}\la \Ze_s -\Ze_{s-},g\ra^2.
\end{equation}
The last sum in the right-hand side is an increasing process with compensator
\begin{equation}\label{mp Xe-comp}
\int_0^t \la \Ze_s,(\sigma^2 +\e \ind_{\G})g^2\ra \ ds.
\end{equation}

\medskip
The proof of Theorem \ref{theo cv to sbm} relies on the following two results, in which we use the notation $(\cY_t)_{t\geq 0}$ for the canonical process on $D_{\cM_f(\R^d)}[0,\infty)$. 
Recall from (\ref{scalingtransfo}) the definition of the process $X^{\varpi,\varepsilon}$
in terms of $Z^{\varpi,\varepsilon}$.

\begin{lemma}\label{lemma: tightness}
For every $\vp\in \bOm$ and $\e\in(0,1)$, let $\Pi^{\vp,\e}$ be the  law of the process $\Xe$ under $\P_{[ \e^{-1}] \delta_0}$. Then,
\begin{description}
\item{\rm (i)} For every $\delta>0$ and $T>0$, there is a compact subset $K_{\delta,T}$ of $\R^d$ such that, for every $\varpi\in\bOm$,
$$\sup_{\e\in(0,1)}\Pi^{\vp,\e}\Big(\sup_{0\leq t\leq T} \cY_t(K_{\delta,T}^c)\Big) < \delta$$
\item{\rm (ii)} For every $g\in \bar{C}^2(\R^d)$ and $\vp\in\bOm$, the collection of the laws of the process $(\la \cY_t,g\ra)_{t\geq 0}$ under $\Pi^{\vp,\e}$, $\e\in(0,1)$, is relatively compact in the space of all probability measures on $D_{\R}[0,\infty)$.
\end{description}
Consequently, for every $\vp\in\bOm$, the collection $(\Pi^{\vp,\e})_{\e\in(0,1)}$ is relatively compact in the space of all probability measures on $D_{\cM_f(\R^d)}[0,\infty)$.
\end{lemma}

The last assertion of the lemma is an immediate consequence of (i) and (ii) using Theorem II.4.1 in \cite{PER2002}.

\begin{prop}\label{prop: conv fd}
Let $g\in\bar{C}_+^2(\R^d)$. There exists a measurable subset $\bOm_{g}$ of $\bOm$ such that $\bP(\bOm_g)=1$ and the following holds for every $\varpi\in \bOm_g$. For every $s,t\geq 0$, for every integer $p\in \N$ and every choice of $t_1,\ldots,t_p\in[0,t]$ and $f_1,\ldots,f_p\in \bar{C}\big(\cM_f(\R^d)\big)$, we have
$$
\lim_{\e\rightarrow 0}\E_{[ \e^{-1}] \delta_0}\bigg[\bigg\{e^{-\la \Xe_{t+s},g\ra}-e^{-\la \Xe_t,g\ra}-\int_t^{t+s}\Big\la \Xe_u,-\frac{1}{2}\Delta g + \psi_{(\kappa)}(g)\Big\ra \ e^{-\la \Xe_u,g\ra}du \bigg\} \prod_{i=1}^p f_i(\Xe_{t_i})\bigg]=0.
$$
\end{prop}

We postpone the proof of Lemma \ref{lemma: tightness} and Proposition \ref{prop: conv fd},
and explain how
Theorem \ref{theo cv to sbm} follows from these two statements. We choose a countable dense subset $G$ of $\bar{C}_+^2(\R^d)$ and set
$$\bOm'=\bigcap_{g\in G} \bOm_g.$$
Fix $\vp\in\bOm'$. By Lemma \ref{lemma: tightness}, the collection $(\Pi^{\vp,\e})_{\e\in(0,1)}$ is relatively compact. Let $\Pi^*$ be a sequential limit of this collection as $\e$ tends to $0$. We deduce from Proposition \ref{prop: conv fd} that, for every $g\in G$, for every $s,t\geq 0$ and every choice of $t_1,\ldots,t_p\in[0,t]$ and $f_1,\ldots,f_p\in \bar{C}\big(\cM_f(\R^d)\big)$, we have 
\begin{equation}
\label{MPlimit}
\Pi^*\bigg[\bigg\{e^{-\la \cY_{t+s},g\ra}-e^{-\la \cY_t,g\ra}-\int_t^{t+s}\Big\la \cY_u,-\frac{1}{2}\Delta g + 
\psi_{(\kappa)}(g)\Big\ra \ e^{-\la \cY_u,g\ra}du \bigg\} \prod_{i=1}^p f_i(\cY_{t_i})\bigg]=0.
\end{equation}
The required passage to the limit under the expectation sign is easily justified: Note that, by comparing with the case when there is no killing and using standard results of the theory of critical branching processes, we have for every $T>0$
\begin{equation}
\label{unibound}
\sup_{\e\in(0,1)} \Big( \E_{[\e^{-1}\delta_0]}\Big[\sup_{0\leq r\leq T} \la X^{\vp,\e}_r,1\ra^2\Big]\Big) <\infty.
\end{equation}
Since $G$ is dense in $\bar{C}_+^2(\R^d)$, another easily justified passage to the limit shows that (\ref{MPlimit}) holds for every $g\in \bar{C}_+^2(\R^d)$. Thus, $\Pi^*$ satisfies the martingale problem for super-Brownian motion with branching mechanism $\psi_{(\kappa)}$ as stated in Section \ref{section: preliminaries}. Since it is also clear that $\Pi^*(Y_0=\delta_0)=1$, $\Pi^*$ must be the law of super-Brownian motion with branching mechanism $\psi_{(\kappa)}$ started from $\delta_0$. This completes the proof of Theorem \ref{theo cv to sbm},
but we still need to establish Lemma \ref{lemma: tightness} and Proposition \ref{prop: conv fd}. $\hfill\square$

\medskip
\noindent\textbf{Proof of Lemma \ref{lemma: tightness}. } The compact containment condition (i) in the lemma is immediately obtained by observing that $\Xe$ is dominated by $X^0=X^{\varpi,0}$ and by using the convergence of rescaled branching Brownian motions (without killing in the obstacles) towards super-Brownian motion. So we only need to verify (ii). In the remaining part of the proof, we fix $\vp\in\bOm$ and $g\in \bar{C}^2_+(\R^d)$. To simplify notation, we set 
$${\cal X}^\e_t=\la \Xe_t,g \ra.$$ 
By the remarks following Lemma \ref{lemma: mart pb} and an elementary scaling transformation, we have
$${\cal X}^\e_t ={\cal X}^\e_0 +  M^\e_t + V^\e_t,$$
where 
$$V^\e_t=\int_0^t \Big\la \Xe_s,\frac{1}{2}\Delta g - \ind_{\Ge}g\Big\ra\ ds$$
and $M^\e$ is a martingale, whose square bracket is given by
$$[M^\e,M^\e]_t= \int_0^t \la \Xe_s,\e \nabla g.\nabla g\ra\ ds +\sum_{0\leq s\leq t}\la \Xe_s -\Xe_{s-},g\ra^2.$$
Furthermore, the oblique bracket of $M^\e$ is equal to
$$\la M^\e,M^\e\ra_t = \int_0^t \big\la \Xe_s,\sigma^2g^2+\e\big(\nabla g.\nabla g+\ind_{\Ge}g^2\big)\big\ra \ ds.$$
By standard criteria (see in particular Theorem VI.4.13 in Jacod and Shiryaev \cite{JS}), the tightness of the laws of the processes ${\cal X}^\e$, $\e\in(0,1)$, will follow if we can verify that the laws of the processes $V^\e$ and $\la M^\e,M^\e\ra$, for $\e\in(0,1)$, are $C$-tight. But this is immediate from the preceding explicit formulas and (\ref{unibound}). $\hfill\square$

\medskip
\noindent\textbf{Proof of Proposition \ref{prop: conv fd}. }Let us fix $s,t\geq 0$, $t_1,\ldots,t_p\in[0,t]$ and $f_1,\ldots,f_p \in \bar{C}\big(\cM_f(\R^d)\big)$. Let $\vp \in \bOm$ and $\e\in (0,1)$. Using the facts that $\Upsilon(1)=\Upsilon'(1)=1$, $\Upsilon''(1)=\sigma^2$ and Lemma \ref{lemma: mart pb} applied to the function $f(x)=\exp\big\{- \e g(x\sqrt{\e})\big\}$, we can write
\begin{eqnarray}
0&=& \Ee\bigg[\bigg\{e^{-\la \Xe_{t+s},g \ra}-e^{-\la \Xe_t,g \ra} -\e^{-2}\int_t^{t+s}\Big\la \Xe_u,\frac{1}{2}\big(- \e^2 \Delta g + \e^3\nabla g.\nabla g\big) \nonumber\\
& & +\ \e \ind_{\Ge} \big(1-e^{-\e g}\big)e^{\e g} + e^{\e g}\big(\Upsilon\big(e^{-\e g}\big)- e^{-\e g}\big)\Big\ra \ e^{-\la \Xe_u, g\ra}du  \bigg\}\prod_{i=1}^pf_i\big(\Xe_{t_i}\big)\bigg] \nonumber\\
&=& \Ee\bigg[\bigg\{e^{-\la \Xe_{t+s},g \ra}-e^{-\la \Xe_t,g \ra} -\int_t^{t+s}\Big\la \Xe_u,-\frac{1}{2} \Delta g +\frac{\e}{2}\nabla g.\nabla g \nonumber \\
& & +\ \ind_{\Ge}(g +\eta_{\e}) +\frac{\sigma^2}{2} (g^2 +\tilde{\eta}_{\e})\Big\ra \ e^{-\la \Xe_u,g\ra}du \bigg\}\prod_{i=1}^pf_i\big(\Xe_{t_i}\big)\bigg], \label{unif conv}
\end{eqnarray}
where $\|\eta_{\e}\|\leq C_1(g) \e $ for some constant $C_1(g)$ depending only on $g$, and $\|\tilde{\eta}_{\e}\|\rightarrow 0$ as $\e\rightarrow 0$. On the one hand, using Fubini's theorem, the fact that $(\la \Xe_{r},1\ra)_{r\geq 0}$ is a supermartingale and the inequality $\Ee[\la\Xe_0,1\ra]\leq 1$, we have \setlength\arraycolsep{1pt}
\begin{eqnarray*}
\bigg|\Ee&\bigg[&\bigg\{\int_t^{t+s}\Big\la \Xe_u,\frac{\e}{2} \nabla g .\nabla g + \ind_{\Ge}\eta_{\e} +\frac{\sigma^2}{2} \tilde{\eta}_{\e}\Big\ra
\ e^{-\la \Xe_u,g\ra}du\bigg\}\prod_{i=1}^p f_i\big(\Xe_{t_i}\big)\bigg]\bigg| \\
&\leq& \bigg(\frac{\e}{2}\big(\|\nabla g.\nabla g\|+2C_1(g)\big)+ \frac{\sigma^2}{2} \|\tilde{\eta}_{\e}\|\bigg) \bigg(\prod_{i=1}^p\|f_i\|\bigg)\Ee\bigg[\int_t^{t+s}\la \Xe_u,1\ra\ du\bigg]\\
&\leq &\big(C_2(g)\,\e + (\sigma^2/2)\|\tilde{\eta}_{\e}\|\big) 
\ \bigg(\prod_{i=1}^p\|f_i\|\bigg)\,s
\end{eqnarray*}
with a constant $C_2(g)$ depending only on $g$. On the other hand, \setlength\arraycolsep{1pt}
\begin{eqnarray}
\Ee&\bigg[&\bigg\{\int_t^{t+s}\big\la \Xe_u,\ \ind_{\Ge} g\big\ra e^{-\la \Xe_u, g\ra}du\bigg\}\prod_{i=1}^p f_i\big(\Xe_{t_i}\big)\bigg]\nonumber \\
&=& \Ee\bigg[\bigg\{\int_t^{t+s}\la \Xe_u,\ka  g\ra\ e^{-\la \Xe_u, g\ra}du\bigg\}\prod_{i=1}^p f_i\big(\Xe_{t_i}\big)\bigg]\nonumber \\
& & +\ \Ee\bigg[\bigg\{\int_t^{t+s}\big\la \Xe_u,\big(\ind_{\Ge}-\ka\big) g\big\ra\ e^{-\la \Xe_u, g\ra}du\bigg\}\prod_{i=1}^p f_i\big(\Xe_{t_i}\big)\bigg].
\label{remainder mp-N}
\end{eqnarray}
The absolute value of the second term in the right-hand side of (\ref{remainder mp-N}) is bounded above by
$$
\bigg(\prod_{i=1}^p\|f_i\|\bigg)\int_t^{t+s}\Ee\Big[\big|\big\la \Xe_u, \big(\ind_{\Ge}-\ka\big)g \big\ra\big|\Big]du.
$$
Consequently, going back to (\ref{unif conv}) and using the preceding estimates, we obtain for every $\vp\in \bOm$ and $\e\in(0,1)$,
\begin{eqnarray}
\label{remainderMP}
&&\bigg| \Ee\bigg[\bigg\{e^{-\la \Xe_{t+s}, g\ra}-e^{-\la \Xe_t, g\ra}-\int_t^{t+s}\Big\la \Xe_u,-\frac{1}{2}\Delta g + \ka g+\frac{\sigma^2}{2} g^2\Big\ra\ e^{-\la \Xe_u, g\ra}du\bigg\}\prod_{i=1}^p f_i\big(\Xe_{t_i}\big)\bigg]\bigg|\nonumber\\
&&\qquad\leq \bigg(\prod_{i=1}^p\|f_i\|\bigg)\bigg(\big\{C_2(g)\,\e + (\sigma^2/2)\|\tilde{\eta}_{\e}\|\big\} s + r_{\e}(\vp,g,t,t+s)\bigg),
\end{eqnarray}
where $r_{\e}(\vp,g,t_1,t_2)$ is defined for $0\leq t_1\leq t_2$ by
$$
r_{\e}(\vp,g,t_1,t_2)=\int_{t_1}^{t_2}\Ee\Big[\big|\big\la \Xe_u, \big(\ind_{\Ge}-\ka\big)g \big\ra\big|\Big]du.
$$

\begin{lemma}\label{lemm mixing} Let $u> 0$. Then, $\bP$-a.s.
\begin{equation}\label{eq mixing}
\lim_{\e\rightarrow \infty} \Ee\Big[\big|\big\la \Xe_u,\ \big(\ind_{\Ge} -\ka \big) g\big\ra \big|\Big]=0.
\end{equation}
\end{lemma}

Assuming that the lemma is proved, we can readily complete the proof of Proposition \ref{prop: conv fd}. Using Fubini's theorem,
we may find a set $\bOm_g$ with $\bP\big[\bOm_g\big]=1$, such that, for every $\varpi\in\bOm_g$, the convergence in (\ref{eq mixing}) holds simultaneously for all $u> 0$, except possibly on a set of zero Lebesgue measure (depending on $\vp$). Since  for every $\vp\in\bOm$, $\e\in (0,1)$ and $u\geq 0$
$$
\Ee\Big[\big|\big\la \Xe_u,\ \big(\ind_{\Ge}-\ka \big)g\big\ra \big|\Big]\leq \| g\|\, \Ee[\la \Xe_u,1\ra]\leq \|g\|,
$$
dominated convergence guarantees that for each $\vp\in \bOm_g$, $\lim_{\e\rightarrow 0}r_{\e}(\vp,g,t,t+s)=0$ for all $t,s\geq 0$.
It follows that the right-hand side of (\ref{remainderMP}) tends to $0$ as $\e\to 0$ when $\varpi\in\bOm_g$, which completes the proof of Proposition \ref{prop: conv fd}. $\hfill\square$

\medskip
\noindent\textbf{Proof of Lemma \ref{lemm mixing}. } Recall that $\xi^{\varpi,\e}$ denotes standard $d$-dimensional Brownian motion killed at rate $\e$ within $\Gamma_\varpi$. We also use the notation $\kie$ for Brownian motion killed at rate $1$ within $\Ge$. Both processes 
$\xi^{\varpi,\e}$ and $\kie$ start from $x$ under the probability measure $\rmP_x$. 

We first recall classical moment formulas for branching Brownian motion. For every $x\in\R^d$,
$k\in\N$, and every bounded measurable function $h$ on $\R^d$, we have
$$\E_{k\delta_x}[\la Z^{\varpi,\e}_t,h\ra] = k\,\rmE_x[h(\xie_t)]$$
and 
$$\E_{k\delta_x}[\la Z^{\varpi,\e}_t,h\ra^2] = k\,\rmE_x[h(\xie_t)^2]
+ k(k-1)\rmE_x[h(\xie_t)]^2 + k\sigma^2\,\rmE_x\left[\int_0^t ds\,
\rmE_{\xie_s}[h(\xie_{t-s})]^2\right].$$
These formulas are easily derived, first for $k=1$, from the well-known formula for the
Laplace functional of $Z^{\varpi,\e}_t$ under $\P_{\delta_x}$: See e.g. Proposition II.3
in \cite{LEG1999}. 

Recalling the definition of $X^{\varpi,\e}$ in terms of $Z^{\varpi,\e}$, we get similar formulas
for the first and second moment of $\la X^{\varpi,\e}_s,h\ra$. In particular, for $s\geq0$
and for every $\vp \in \bOm$ and $\e\in (0,1)$, we have
\begin{equation}
\Ee\Big[\big\la \Xe_s,\big(\ind_{\Ge}-\ka\big)g \big\ra^2\Big]= \e^2\big([ \e^{-1}]^2-[ \e^{-1}]\big) A^{\vp,\e}_1(0,s,g)^2
+ \e[ \e^{-1}] A^{\vp,\e}_2(0,s,g), \label{log-laplace bbm}
\end{equation}
where, for every $x\in\R^d$,
$$
A^{\vp,\e}_1(x,s,g):= \rmE_x\Big[\big(\ind_{\Ge}(\kie_s)-\ka\big)g(\kie_s)\Big], $$
and 
$$A^{\vp,\e}_2(x,s,g):= \sigma^2\ \rmE_x\bigg[\int_0^s A^{\vp,\e}_1(\kie_v,s-v,g)^2 dv\bigg] + \e\ \rmE_x\Big[\big(\ind_{\Ge}(\kie_s)-\ka\big)^2g(\kie_s)^2\Big].
$$

We claim that, for every $x\in \R^d$ and $s> 0$, 
\begin{equation}\label{conv bm}
\bP\mathrm{-a.s.},\qquad \lim_{\e \rightarrow 0}\rmE_x\Big[\big(\ind_{\Ge}\big(\kie_s\big)-\ka\big)g\big(\kie_s\big)\Big]=0.
\end{equation}

Assume for the moment that the claim holds. By taking $x=0$ in (\ref{conv bm}), we obtain that $A^{\vp,\e}_1(0,s,g)$ tends to $0$ as $\e\rightarrow 0$, $\bP$-a.s. 
Consider next $A^{\vp,\e}_2(0,s,g)$. The second term in the formula for $A^{\vp,\e}_2(x,s,g)$ obviously tends to $0$ uniformly in $x,s$ and independently of the obstacles, as $\e\rightarrow 0$. To handle the first term, use Fubini's theorem to obtain that, for every $\varpi$ belonging to a set $\widetilde\bOm$ of full probability, there exists a set $\cN_{\vp}\subset \R^d\times \R_+$ of zero Lebesgue measure such that the convergence in (\ref{conv bm}) holds simultaneously for all $(x,s)\notin \cN_{\vp}$. We can then write the first term in the formula for $A^{\vp,\e}_2(0,s,g)$ as follows: 
\begin{equation}\label{branching term}
\sigma^2 \int_0^s dv\int_{\R^d}dy\ p_v^{\vp,\e}(0,y)\ \rmE_y\Big[ \big(\ind_{\Ge}(\kie_{s-v})-\ka\big)g(\kie_{s-v})\Big]^2,
\end{equation}
where $p^{\vp,\e}_v(\cdot,\cdot)$ is the transition kernel of $\kie$ at time $v$. Plainly, we have $p^{\vp,\e}_v(x,y)\leq p_v(x,y)$ for each $v\geq 0$ and $\e \in (0,1)$, where $p_v(\cdot,\cdot)$ is the transition kernel of standard Brownian motion at time $v$. Using this bound and dominated convergence gives us that the quantity in (\ref{branching term}) tends to $0$ as $\e\to 0$, for every $\varpi\in\widetilde \bOm$. Hence $A^{\vp,\e}_2(0,s,g)$ tends to $0$ as $\e\to 0$, $\bP$-a.s. The result of Lemma \ref{lemm mixing} now follows from (\ref{log-laplace bbm}) and the Cauchy-Schwarz inequality.

We still have to prove our claim (\ref{conv bm}). We fix $x\in \R^d$ and $s> 0$. Let $\vp\in \bOm$ and $\e \in (0,1)$. From the definition of Brownian motion killed at rate $1$ in $\Ge$, we have
\begin{equation}\label{FK}
\rmE_x\Big[\big(\ind_{\Ge}(\kie_s)-\ka\big)g(\kie_s)\Big] =\rmE_x\Big[\big(\ind_{\Ge}(\xi_s)-\ka\big)g(\xi_s)\exp\Big\{-\int_0^s\ind_{\Ge}(\xi_u)du\Big\} \Big].
\end{equation}
Let $\eta>0$ and choose $\theta\in(0,s)$ such that $e^{-\theta}\geq 1-\eta$. By the Markov property applied at time $s-\theta$, the quantities in (\ref{FK}) are equal to \setlength\arraycolsep{1pt}
\begin{eqnarray}
\rmE_x&\bigg[&\exp\Big\{-\int_0^{s-\theta}\ind_{\Ge}(\xi_u)du\Big\}\ \rmE_{\xi_{s-\theta}}\Big[\big(\ind_{\Ge}(\xi_{\theta})-\ka\big)g(\xi_{\theta}) \exp\Big\{- \int_0^{\theta}\ind_{\Ge}(\xi_v)dv\Big\}\Big] \bigg] \nonumber\\
&=& \rmE_x\bigg[\exp\Big\{-\int_0^{s-\theta}\ind_{\Ge}(\xi_u)du\Big\}\ \rmE_{\xi_{s-\theta}}\Big[\big(\ind_{\Ge}(\xi_{\theta})-\ka\big) g(\xi_{\theta})\Big] \bigg]\label{decorrel} \\
& & + \rmE_x\bigg[\exp \Big\{-\int_0^{s-\theta}\ind_{\Ge}(\xi_u)du\Big\}\ \rmE_{\xi_{s-\theta}} \Big[\big(\ind_{\Ge}(\xi_{\theta})-\ka\big)g(\xi_{\theta})\Big(\exp
\Big\{-\int_0^{\theta}\ind_{\Ge}(\xi_v)dv\Big\}-1\Big)\Big] \bigg].
\nonumber
\end{eqnarray}
The condition $1-e^{-\theta} \leq \eta$ entails that the absolute value of the second term in the right-hand side of (\ref{decorrel}) is bounded above by $\|g\|\eta$. Suppose we know that for every $y\in\R^d$ we have
\begin{equation}\label{key conv}
\bP\mathrm{-a.s.},\qquad \lim_{\e \rightarrow 0}\rmE_y\Big[\big(\ind_{\Ge}(\xi_{\theta})-\ka\big)g(\xi_{\theta})\Big]=0.
\end{equation}
Then, using the fact that the law of $\xi_{s-\theta}$ is absolutely continuous with respect to Lebesgue measure
(and Fubini's theorem to see that the convergence in (\ref{key conv}) holds simultaneously for almost all $y\in\R^d$, $\bP$-a.s.), we conclude that the 
first term in the right-hand side of (\ref{decorrel}) converges to $0$ with $\bP$-probability $1$ as $\e\rightarrow 0$. Consequently, we have $\bP$-a.s.
$$
\limsup_{\e\rightarrow 0} \Big|\rmE_x\Big[\big(\ind_{\Ge}(\kie_s)-\ka\big)g(\kie_s)\Big]
\Big|\leq \|g\|\eta.
$$
Since $\eta$ was arbitrary, our claim (\ref{conv bm}) follows.

It thus remains to prove that (\ref{key conv}) holds. We fix $y\in\R^d$ and $\theta>0$. In the following, $\xi'$ stands for another Brownian motion independent of $\xi$, which also starts from $y$ under the probability measure $\rmP_y$. For each $\e\in (0,1)$, we have
\begin{eqnarray}
\bE\Big[\rmE_y\big[\big(\ind_{\Ge}(\xi_{\theta})-\ka\big)g(\xi_{\theta}) \big]^2\Big]&=&\bE\Big[\rmE_y\big[\big(\ind_{\Ge}(\xi_{\theta}) -\ka\big)\big(\ind_{\Ge}(\xi'_{\theta})-\ka\big)g(\xi_{\theta})g(\xi'_{\theta})\big]\Big] \nonumber\\
&=& \rmE_y\Big[g(\xi_{\theta})g(\xi'_{\theta})\Big\{\bP\big[\xi_{\theta}\in \Ge;\ \xi'_{\theta}\in \Ge\big]-\ka^2\Big\}\Big],\label{spatial correl}
\end{eqnarray}
where the last line uses Fubini's theorem and the definition of $\ka$. Recall that the measure $\Theta$ is supported on compact sets which are contained in the fixed ball  $\ov{B}(0,r_0)$. If $|x-x'|>2r_0$, the sets $\{(z,K):x\in z+ K\hbox{ and }K\subset \ov B(0,r_0)\}$
and  $\{(z,K):x'\in z+ K\hbox{ and }K\subset \ov B(0,r_0)\}$ are disjoint, and so the events $\{x\in\G\}$
and $\{x'\in \G\}$ are independent under $\bP$. Recalling that $\Ge=\sqrt{\e}\ \G$, we see that if $|x-x'|>2r_0\sqrt{\e}$ we have $\bP[x\in \Ge;\ x'\in \Ge]=\bP[x\in \Ge]\bP[x'\in \Ge]=\ka^2$, which enables us to write
$$
\bP\big[\xi_{\theta}\in \Ge;\ \xi'_{\theta}\in \Ge\big]=\ka^2\ \ind_{\{|\xi_{\theta}-\xi'_{\theta}|> 2r_0\sqrt{\e}\}}+ \bP\big[\xi_{\theta}\in \Ge;\ \xi'_{\theta}\in \Ge\big]\ \ind_{\{|\xi_{\theta}-\xi'_{\theta}|\leq 2r_0\sqrt{\e}\}}.
$$
Going back to (\ref{spatial correl}), we obtain
\begin{eqnarray}
\bE\Big[\rmE_y\big[\big(\ind_{\Ge}(\xi_{\theta})-\ka\big)g(\xi_{\theta})\big]^2\Big]&=& \rmE_y\Big[g(\xi_{\theta})g(\xi'_{\theta})\ind_{\{|\xi_{\theta}-\xi'_{\theta}|\leq 2r_0\sqrt{\e}\}} \Big\{\bP\big[\xi_{\theta}\in \Ge;\ \xi'_{\theta}\in \Ge\big]-\ka^2\Big\}\Big]\nonumber\\
&\leq & \|g\|^2 \rmP_y\big[|\xi_{\theta}-\xi'_{\theta}|\leq 2r_0\sqrt{\e}\big]\nonumber\\
&=& \|g\|^2 \rmP_0\big[|\xi_{2\theta}|\leq 2r_0\sqrt{\e} \big] \nonumber \\
& \leq &C\frac{\|g\|^2}{\theta^{d/2}}\ \e^{d/2},
\label{expected mixing}
\end{eqnarray}
where the constant $C>0$ depends only on $r_0$. Let $\eta>0$. By the Markov inequality, we can write
$$
\bP\Big[\big|\rmE_y\big[\big(\ind_{\Ge}(\xi_{\theta})-\ka\big)g(\xi_{\theta})\big]\big|>\eta\Big]\leq \frac{1}{\eta^2}\ \bE\Big[\rmE_y\big[\big(\ind_{\Ge}(\xi_{\theta})-\ka\big)g(\xi_{\theta})\big]^2\Big]\leq C\frac{\|g\|^2}{\eta^2\theta^{d/2}}\ \e^{d/2}.
$$
Applying the last bound with $\e=\e_k=k^{-3}$, for every $k\in \N$, yields a convergent series, and so by the Borel-Cantelli lemma we have
$\bP$-a.s.,
$$
\limsup_{k\rightarrow \infty}\ \Big|\rmE_y\big[\big(\ind_{\Gamma_{\vp}^{\e_k}}(\xi_{\theta})-\ka\big)g(\xi_{\theta})\big]\Big|\leq \eta.
$$
Since $\eta>0$ was arbitrary, the convergence (\ref{key conv}) holds along the sequence $(\e_k)_{k\geq 1}$.

To complete the proof, we set for every $\vp\in \bOm$ and $\e\in (0,1)$
$$
\Ue=\rmE_y\big[\ind_{\Ge}(\xi_{\theta})\,g(\xi_{\theta})\big]
$$
(recall that $y$ and $\theta$ are fixed).
We shall prove the following result: For every $\vp\in \bOm$,
\begin{equation}\label{uniformity}
\lim_{k\rightarrow \infty}\sup_{(k+1)^{-3}\leq \e \leq k^{-3}}\big|\Ue-U_{\vp,\e_k}\big|=0.
\end{equation}
Combining (\ref{uniformity}) with the fact that the convergence (\ref{key conv}) 
is true along the sequence $(\e_k)_{k\geq 1}$ will then lead to the desired conclusion. 
In order to prove (\ref{uniformity}), we first note that for every $\e\in(0,1)$,
$$
\Ue=\frac{1}{(2\pi \theta)^{d/2}}\int_{\R^d}\ind_{\Ge}(x)\,g(x)e^{-\frac{|x-y|^2}{2\theta^2}}dx = \int_{\R^d}\ind_{\Ge}(x)\, h(x)dx,
$$
with a function $h\in \bar{C}^2_+(\R^d)$ which depends on $y$ and $\theta$, but not on $\vp$. Furthermore, for any fixed $\eta>0$ we can find a large closed ball $B$ centered at the origin and such that $\int_{B^c}h(x)dx<\eta$. Hence, if we set
$$
\Ue'=\int_B h(x)\ \ind_{\Ge}(x)dx,
$$
we have $|\Ue-\Ue'|\leq \eta$ for every $\e\in(0,1)$. Thanks to this remark, it is enough to prove that (\ref{uniformity}) holds when $\Ue$ and $U_{\vp,\e_k}$ are replaced
 by $\Ue'$ and $U'_{\vp,\e_k}$ respectively.

Let $k\in \N$ and $\e\in [(k+1)^{-3},k^{-3})$. We have
\begin{eqnarray}
\Ue'&=& \e^{d/2}\int_{\e^{-1/2}B}h\big(x\sqrt{\e}\big)\ind_{\G}(x)dx \nonumber\\
&=& \e^{d/2}\int_{k^{3/2}B}h\big(x \sqrt{\e}\big)\ind_{\G}(x)dx + \e^{d/2}\int_{(\e^{-1/2}B)\setminus
(k^{3/2}B)}h\big(x\sqrt{\e}\big)\ind_{\G}(x)dx \label{eq unif}.
\end{eqnarray}
Now, the first term in the right-hand side of (\ref{eq unif}) is equal to
\setlength\arraycolsep{1pt}
\begin{eqnarray*}
&& \big(\e k^3\big)^{d/2}\frac{1}{k^{3d/2}} \int_{k^{3/2}B}h\Big(\frac{x}{k^{3/2}}\Big)\ind_{\G}(x)dx + \e^{d/2} \int_{k^{3/2}B}\Big[h\Big(x\sqrt{\e}\Big)-h\Big(\frac{x}{k^{3/2}}\Big)\Big]\ind_{\G}(x)dx\\
&&\qquad= \big(\e k^3\big)^{d/2}U'_{\vp,\e_k} + \iota(\vp,\e,k) ,
\end{eqnarray*}
where
$$\iota(\vp,\e,k) = (\e k^3)^{d/2} \int_B \Big[h((k^{3/2}\sqrt{\e})y)-h(y)\Big]\,\ind_{\G}(k^{3/2}y)dy.$$
Note that $0\leq 1-k^{3/2}\sqrt{\e}\leq \frac{C}{k}$ with a constant $C$ independent of $\e$ and $k$, from which it easily follows that
$$
\sup_{(k+1)^{-3}\leq \e \leq k^{-3}}|\iota(\vp,\e,k)|\leq \frac{C'}{k}\ \|\nabla h\|,
$$
with a constant $C'$ depending only on $B$. 

Similarly, the second term in the right-hand side of (\ref{eq unif}) is bounded above by
$$
\e^{d/2}\ \|h\|\big[\e^{-d/2}-k^{3d/2}\big]\lambda_d(B)\leq \frac{C''\|h\|}{k}.
$$
Finally, from (\ref{eq unif}) and the preceding estimates, 
we have
$$\big|\Ue'-U'_{\vp,\e_k}\big|\leq (1-(\e k^3)^{d/2})\,U'_{\vp,\e_k}
+ \frac{C' \|\nabla h\|}{k} + \frac{C''\|h\|}{k},$$
and the convergence (\ref{uniformity}) follows. 
This completes the proof.
$\hfill\square$

\medskip
\noindent{\bf Proof of Proposition \ref{large-dev}}. Clearly we get a lower bound for the quantity $\P_{\delta_0}(Z^{\vp,\e} \hbox{ hits } B(0,R)^c)$ if we replace $Z^{\vp,\e}$
by branching Brownian motion killed homogeneously over $\R^d$ at rate $\e$. Part (i) of the proposition thus follows 
from the lower bound in Lemma \ref{Sawyer}. 

Let us turn to the proof of the first assertion in (ii). We have the following inequality: For every $\vp\in \bOm$ and $\e \in (0,1)$,
\begin{equation}\label{ineq hitting}
\P_{\delta_0}\big[\Ze\ \mathrm{hits}\ B(0, R)^c\big]\leq \rmP_0\big[\xie\ \mathrm{hits}\ B(0,R)^c\big].
\end{equation}
Indeed, using the formalism of Subsection \ref{subs: BBM et ob}, the criticality of the offspring distribution can be used to check that the right-hand side of (\ref{ineq hitting}) is
just the expected value of the number of those historical paths $\omega^v$ that first exit $B(0,R)$ during the interval
$[\alpha_v,\beta_v \wedge \zeta^{\varpi,\e}_v)$. Alternatively, it is easy to derive an integral
equation similar to (\ref{intequ}) for the function $x\to \P_{\delta_x}\big[\Ze\ \mathrm{hits}\ B(0, R)^c\big]$,
and the bound (\ref{ineq hitting}) then trivially follows from this integral equation.

Then, let us bound $\bP\otimes\rmP_0\big[\xie\ \mathrm{hits}\ B(0, R)^c\big]$. 
We write $\xi^{\varpi,\e,k}$ for the $k$-th coordinate of $\xie$, for every $k=1,\ldots,d$. First, observe that 
\begin{equation}\label{coordinates}
\bP\otimes\rmP_0\big[\xie\ \mathrm{hits}\ B(0, R)^c\big] \leq \sum_{k=1}^d \bP\otimes\rmP_0\big[\xi^{\varpi,\e,k}\ \mathrm{hits}\ (- R/\sqrt{d}, R/\sqrt{d})^c\big].
\end{equation}
Clearly, we can restrict our attention to the first term in the sum. Define $N_{\e,R}:=[ R\sqrt{\e}/\sqrt{d}]$ and for every $j\in \Z$, set $T_j^{\e}
=\inf\{t\geq 0:\xi^1_t=j\e^{-1/2}\}$, where $\xi^1$ stands for the first coordinate of $\xi$. In this notation, we can write
\begin{eqnarray}
&\bP&\otimes\  \rmP_0\big[\xi^{\varpi,\e,1}\ \mathrm{hits}\ (- R/\sqrt{d}, R/\sqrt{d})^c\big]\nonumber\\
&\leq & \bP\otimes\rmP_0\Big[\xi^{\varpi,\e,1}\ \mathrm{hits}\ \big(-N_{\e,R}\e^{-1/2},N_{\e,R}\e^{-1/2}\big)^c\Big]\nonumber \\
&\leq & \bE\otimes \rmE_0\bigg[\exp\bigg\{-\e\int_0^{T_{N_{\e,R}}^{\e}}\ind_{\G}(\xi_s)ds\bigg\}\bigg]+ \bE\otimes \rmE_0\bigg[\exp\bigg\{-\e\int_0^{T_{-N_{\e,R}}^{\e}}\ind_{\G}(\xi_s)ds\bigg\}\bigg]. \label{decoupage ld}
\end{eqnarray}
Consider the first term in the right-hand side of (\ref{decoupage ld}). 
Recall the definition (\ref{obstacleset}) of the set of obstacles, and
set for every $j\in\N$ and $\e>0$,
$$
I_j^{\e}:=\big\{i\in I:\ x_i\in ((j-1)\e^{-1/2},j\e^{-1/2})\times \R^{d-1}\big\},
$$
$$
\G(j,\e):=\bigcup_{i\in I_j^{\e}} (x_i+K_i)\qquad \mathrm{and} \qquad \Ge(j):=\sqrt{\e}\ \G(j,\e).
$$
Note that the random sets $\G(j,\e)$, $j\in \N$, are independent under $\bP$, by properties of Poisson measures.
We have then
\begin{eqnarray}
\bE\otimes \rmE_0\bigg[\exp\bigg\{-\e\int_0^{T_{N_{\e,R}}^{\e}}\ind_{\G}(\xi_s)ds\bigg\}\bigg]
&\leq & \bE\otimes \rmE_0\bigg[\exp\bigg\{-\e\sum_{j=1}^{N_{\e,R}}\int_{T_{j-1}^{\e}}^{T_j^{\e}}\ind_{\G(j,\e)}(\xi_s)ds\bigg\}\bigg]\nonumber \\
&=& \bE\otimes \rmE_0\bigg[\exp\bigg\{-\e\int_0^{T_1^{\e}}\ind_{\G(1,\e)}(\xi_s)ds\bigg\}\bigg]^{N_{\e,R}},\label{def alpha}
\end{eqnarray}
where the equality comes from an application of the strong Markov property of $\xi$, together with the independence
of the random sets  $\G(j,\e)$ and the fact that the distribution of each of these random sets is invariant under translations
by elements of $\{0\}\times \R^{d-1}$. By scaling, if $T_1$ denotes the entrance time of $\xi$ into $[1,\infty)\times \R^{d-1}$, we can write
$$
\alpha_{\e}:= \bE\otimes \rmE_0\bigg[\exp\bigg\{-\e\int_0^{T_1^{\e}}\ind_{\G(1,\e)}(\xi_s)ds\bigg\}\bigg]=\bE\otimes \rmE_0\bigg[\exp\bigg\{-\int_0^{T_1}\ind_{\Ge(1)}(\xi_s)ds\bigg\}\bigg].
$$
We then observe that
\begin{equation}
\label{convdist}
\int_0^{T_1}\ind_{\Ge(1)}(\xi_s)\ind_{(0,1)}(\xi^1_s)\,ds \ \build{\longrightarrow}_{\e\to 0}^{(P)}\ \kappa \int_0^{T_1} \ind_{(0,1)}(\xi^1_s)ds
\end{equation}
where the notation $\build{\longrightarrow}_{}^{(P)}$ refers to convergence in probability under $\bP\otimes \rmP_0$. To see
this, we use arguments similar to the proof of Lemma \ref{lemm mixing}. 
Notice that $\bP[y\in \Ge(1)]\leq \kappa$ if $y\in (0,1)\times \R^{d-1}$, with equality if $y\in (r_0\sqrt{\e},1-r_0\sqrt{\e})\times \R^{d-1}$.
Using this remark, and the same argument as in the derivation of (\ref{expected mixing}), we can write for every fixed $u>0$,
\begin{eqnarray*}
\bE &\otimes &\rmE_0\bigg[\bigg\{\int_0^{T_1\wedge u}\big(\ind_{\Ge(1)}(\xi_s)-\ka\big)\ind_{(0,1)}(\xi_s^1)ds\bigg\}^2\bigg]\\
& =& \bE\otimes \rmE_0\bigg[\int_0^{T_1\wedge u}\!\!\int_0^{T_1\wedge u}\big(\ind_{\Ge(1)}(\xi_s)-\kappa)(\ind_{\Ge(1)}(\xi_t)-\ka)
\,\ind_{(0,1)}(\xi_s^1)\,\ind_{(0,1)}(\xi_t^1)\,dsdt\bigg] \\
&=&  \rmE_0\bigg[\int_0^{T_1\wedge u}\!\!\int_0^{T_1\wedge u}\big\{\bP\big[\xi_s,\xi_t\in \Ge(1)\big]-\ka^2\big\}\,\ind_{(0,1)}(\xi_s^1)\,\ind_{(0,1)}(\xi_t^1)\,dsdt\bigg]+O(\e^{1/2}) \\
&=& \rmE_0\bigg[\int_0^{T_1\wedge u}\!\!\int_0^{T_1\wedge u}\big\{\bP\big[\xi_s,\xi_t\in \Ge(1)\big]-\ka^2\big\}\,\ind_{\{|\xi_s-\xi_t|\leq 2r_0\sqrt{\e}\}}\,
\ind_{(0,1)}(\xi_s^1)\,\ind_{(0,1)}(\xi_t^1)\,dsdt\bigg] +O(\e^{1/2})\\
&\leq & \rmE_0\bigg[\int_0^{u}\!\!\int_0^{u}\ind_{\{|\xi_s-\xi_t|\leq 2r_0\sqrt{\e}\}}\,dsdt\bigg] +O(\e^{1/2})
\end{eqnarray*}
where the error term $O(\e^{1/2})$ corresponds to the contribution of times $s,t$ such that $\xi^1_s$ or $\xi^1_t$ belongs to the set
$(0,r_0\sqrt{\e}]\cup[1-r_0\sqrt{\e},1)$. The preceding quantity tends to $0$ as $\e\rightarrow 0$, which yields the convergence (\ref{convdist}). Since
the limiting variable in (\ref{convdist}) is (strictly) positive a.s., we can find $\e_1\in(0,1)$ and $c_1<1$ such that $\alpha_\e\leq c_1$ for every $\e\in(0,\e_1)$. 
Using (\ref{coordinates}), (\ref{decoupage ld}) and (\ref{def alpha}), we arrive at
$$
\bP\otimes \rmP_0\big[\xie\ \mathrm{hits}\ B(0, R)^c\big] \leq 2d\,c_1^{N_{\e,R}},
$$
for $\e\in(0,\e_1)$. This completes the proof of the first assertion in (ii). 

Let us turn to the second assertion. Let $C_2>0$ be a positive constant whose choice will
be specified later. By simple comparison 
arguments, it is enough to prove the desired estimate when $R$
is of the form
$R=2^k$, for $k\in \N$ large enough, and $\e$ is of the form $\e=2^{-j}$, with $j\in\{0,1,\ldots,2k\}$ such that $\e \geq C_2(\log\log R)^2/R^2$. 

By the first assertion in (ii)
and the Markov inequality, 
\begin{equation}
\label{quenchedtech}
\bP\Big[ \P_{\delta_0}\big[\Ze\ \mathrm{hits}\ B(0, R)^c\big]
\geq \exp(-C_1 R\sqrt{\e}/2)\Big] \leq \exp(-C_1 R\sqrt{\e}/2).
\end{equation}
However, if $R=2^k$ and $\e \geq C_2(\log\log R)^2/R^2$, we have
$$R\sqrt{\e} \geq \sqrt{C_2}\,\log\log R= \sqrt{C_2} (\log k + \log\log 2).$$
Using this bound, we can choose the constant $C_2$ sufficiently large so that we get 
a convergence series when we sum the right-hand side of (\ref{quenchedtech})
over all $R=2^k$ and $\e=2^{-j}$ for $j\in\{0,1,\ldots,2k\}$ such that $\e \geq C_2(\log\log R)^2/R^2$.
The Borel-Cantelli lemma now yields the desired result.

$\hfill\square$

\section{Proof of the main result}\label{section: main result}

In this section, we prove Theorem \ref{theo charging proba}. We fix the environment $\varpi$ such that the weak convergence of Theorem \ref{theo cv to sbm} holds, and derive the convergence in Theorem \ref{theo charging proba} for this fixed value of the environment. For the sake of simplicity, we shall omit $\varpi$ in the notation and write $Z^\e$ instead of $Z^{\varpi,\e}$, and $X^\e$ instead of $X^{\varpi,\e}$. 

We shall verify that for any increasing sequence $(R_n)_{n\geq 1}$ of positive reals converging to $+\infty$ and any sequence $(\e_n)_{n\geq 1}$ of nonnegative reals such that $\e_nR_n^2 \rightarrow a \in[0,\infty]$, we have
\begin{equation}
\label{keylimit}
\lim_{n\to\infty} R_n^2\, \P_{\delta_0}(Z^{\e_n}\hbox{ hits }R_nA^c)
=u_{(\kappa a)}(0),
\end{equation}
where $u_{(\infty)}(0)=0$ by convention.

The statement of Theorem  \ref{theo charging proba} follows from this convergence. Indeed, if the conclusion of the theorem fails, then we can find a sequence $R_n\uparrow\infty$ and a sequence $(\e_n)$ of nonnegative reals such that, for every $n\geq 1$, 
$$|R_n^2\, \P_{\delta_0}(Z^{\e_n}\hbox{ hits }R_nA^c) - u _{(\kappa \e_n R_n^2)}(0)| \geq \delta$$
for some constant $\delta >0$. By extracting a subsequence, we may assume that $\e_n R_n^2\longrightarrow a\in[0,\infty]$ and thus obtain a contradiction with
(\ref{keylimit}) since we know from Lemma \ref{analytic} that the mapping $b\to u_{(b)}(0)$ is continuous on $[0,\infty]$.

In proving (\ref{keylimit}), we may assume that $\e_n\to 0$ as $n\to\infty$. Indeed, suppose that (\ref{keylimit}) holds in this particular case and let $(\e'_n)$ be a sequence that does not converge to $0$. If the sequence $\e'_nR_n^2$ converges then necessarily its limit is $+\infty$, and we can find another sequence $\e''_n$ such that $0\leq \e''_n\leq \e'_n$, $\e''_n\to 0$ and $\e''_nR_n^2 \to \infty$. So, if we know that (\ref{keylimit}) holds in the case when the sequence $(\e_n)$ tends to $0$, we obtain 
$$\lim_{n\to\infty}R_n^2\, \P_{\delta_0}(Z^{\e''_n}\hbox{ hits }R_nA^c)=0.$$
However, from the inequality $\e''_n\leq \e'_n$ and a coupling argument (obvious if one uses the construction described in Subsection \ref{subs: BBM et ob}), we get the same result for the sequence $(\e'_n)$. 

A similar comparison argument shows that it is enough to prove (\ref{keylimit}) in the case when $a<\infty$. Otherwise, it suffices to replace $\e_n$ by $\e_n\wedge b R_n^{-2}$ and let $b\to\infty$, using the fact that $u_{(b)}(0)\to 0$ as $b\to\infty$.

Let us now proceed to the proof of (\ref{keylimit}). We fix the sequences $R_n\uparrow \infty$ and $\e_n\to 0$ such that $\e_n R_n^2 \longrightarrow a \in [0,\infty)$. We first assume that $a>0$. The case $a=0$ will be discussed at the end of the section.

Let $B$ be a closed subset of $\R^d$. For every $\e>0$, we have by the definition of $X^\e$
\begin{eqnarray*}
\P_{[\e^{-1}]\delta_0}(X^\e\hbox{ hits } B)
&=& \P_{[\e^{-1}]\delta_0}(\exists t\geq 0 : X^\e_t(B)>0)\\
&=& \P_{[\e^{-1}]\delta_0}\bigg(\exists t\geq 0 : \int Z^\e_{\e^{-1}t}(dx)\,\ind_B(x\sqrt{\e}) >0 \bigg)\\
&=& \P_{[\e^{-1}]\delta_0}(Z^\e \hbox{ hits } \e^{-1/2}B)\\
&=& 1 - \P_{\delta_0}(Z^\e \hbox{ does not hit } \e^{-1/2}B)^{[\e^{-1}]},
\end{eqnarray*}
since (for a fixed environment) the law of $Z^\e$ under $\P_{[\e^{-1}]\delta_0}$ is obtained by adding $[\e^{-1}]$ independent copies of $Z^\e$ under $\P_{\delta_0}$. Applying the preceding identity with $\e=\e_n$ and $B=b_nA^c$, where $b_n=\e_n^{1/2}R_n$, gives us that
\begin{equation}
\label{keyltech1}
1-\P_{\delta_0}(Z^{\e_n}\hbox{ does not hit }R_nA^c)^{[\e_n^{-1}]} =  \P_{[\e_n^{-1}]\delta_0}(X^{\e_n}\hbox{ hits }b_nA^c)
\end{equation}
By Theorem \ref{theo cv to sbm}, we know that the law of $X^{\e_n}$ under $\P_{[\e_n^{-1}]\delta_0}$ converges as $n\to\infty$ to the law of super-Brownian motion with branching mechanism $\psi_{(\kappa)}$ started from $\delta_0$. The next lemma is essentially a consequence of this convergence. We use the notation of Subsection \ref{subs: hitting prob}.

\begin{lemma}
\label{convhittingpro}
We have 
$$\lim_{n\to\infty}
\P_{[\e_n^{-1}]\delta_0}(X^{\e_n}{\rm\ hits\ }b_nA^c)= P_{\delta_0}(Y^{(\kappa)}{\rm\ hits\ }bA^c),$$
where $b=\sqrt{a}=\lim b_n$.
\end{lemma}

We postpone the proof of Lemma \ref{convhittingpro} and proceed to the proof of (\ref{keylimit}), in the case when $a>0$. By the results recalled in Subsection \ref{subs: hitting prob}, we know that
\begin{equation}
\label{hittingsuperBM}
P_{\delta_0}(Y^{(\kappa)}{\rm\ hits\ }bA^c) = 1 - \exp(-v(0))\;,
\end{equation}
where the function $(v(x),x\in bA)$ is the unique nonnegative solution of the
singular boundary value problem
$$\left\{\begin{array}{ll}
\frac{1}{2}\Delta v = \psi_{(\kappa)}(u)\qquad &\hbox{in }bA\;, \vspace{4pt}\\
u_{|\partial (bA)} = +\infty\;.&
\end{array}
\right.$$
It is immediate to verify that $u_{(\kappa a)}(x)= a\,v(bx)$ for every $x\in A$, and in particular $u_{(\kappa a)}(0)=a\,v(0)$. 

From (\ref{keyltech1}), (\ref{hittingsuperBM}) and Lemma \ref{convhittingpro}, we obtain
$$\lim_{n\to\infty} (1-\P_{\delta_0}(Z^{\e_n}\hbox{ hits }R_nA^c))^{[\e_n^{-1}]}=\exp(-v(0))$$
and thus 
$$\lim_{n\to\infty} \e_n^{-1}\,\P_{\delta_0}(Z^{\e_n}\hbox{ hits }R_nA^c)=v(0),$$
or equivalently, since $\e_nR_n^2 \longrightarrow a$, 
$$\lim_{n\to\infty} R_n^2\,\P_{\delta_0}(Z^{\e_n}\hbox{ hits }R_nA^c)=av(0)=u_{(\kappa a)}(0).$$
This completes the proof of (\ref{keylimit}), in the case $a>0$. $\hfill\square$

\smallskip
\noindent{\bf Proof of Lemma \ref{convhittingpro}.} By replacing $A$ with $bA$, we may and shall assume in this proof that $b=1$. We thus have $b_n\longrightarrow 1$ as $n\to\infty$. We first prove that 
\begin{equation}
\label{lowerboundhitting}
\liminf_{n\to\infty}
\P_{[\e_n^{-1}]\delta_0}(X^{\e_n}{\rm\ hits\ }b_nA^c)\geq  P_{\delta_0}(Y^{(\kappa)}{\rm\ hits\ }A^c).
\end{equation}
By Lemma \ref{hittinglemma}, the events $\{Y^{(\kappa)}\hbox{ hits }A^c\}$ and $\{Y^{(\kappa)}\hbox{ hits } \bar A^c\}$ coincide a.s. We can then find a countable  collection $(\vf_i)_{i\geq 1}$ of continuous functions with compact support contained in $(\bar{A})^c$, such that 
$$\{Y^{(\kappa)}\hbox{ hits }(\bar{A})^c\}=\Big\{\sup_{i\geq 1}\Big(\sup_{t> 0}\ \langle Y^{(\kappa)}_t,\vf_i\rangle\Big) >0\Big\},\qquad P_{\delta_0}\ {\rm a.s.}.$$
Hence, if $(t_j)_{j\geq 1}$ is a sequence dense in $[0,\infty)$, we have
\begin{equation}
\label{hittingtech1}
P_{\delta_0}(Y^{(\kappa)}\hbox{ hits }(\bar{A})^c)=\lim_{N\to\infty}\uparrow P_{\delta_0}\Big(\sup_{1\leq i\leq N}\Big(\sup_{1\leq j\leq N} \langle Y^{(\kappa)}_{t_j},\vf_i\rangle \Big) >0\Big).
\end{equation}

However, Theorem \ref{theo cv to sbm} implies that, for every $N\geq 1$,
\begin{equation}
\label{hittingtech2}
\liminf_{n\to\infty}\ \P_{[\e_n^{-1}]\delta_0}\Big(\sup_{1\leq i\leq N}\Big(\sup_{1\leq j\leq N} \langle X^{\e_n}_{t_j},\vf_i\rangle \Big) >0\Big)
\geq P_{\delta_0}\Big(\sup_{1\leq i\leq N}\Big(\sup_{1\leq j\leq N} \langle Y^{(\kappa)}_{t_j},\vf_i\rangle \Big) >0\Big).
\end{equation}
Recall that $b_n\to 1$, and note that the support of each function $\vf_i$ is at a strictly positive distance of the set $A$. As a consequence, for every fixed $N$, the support of $\vf_i$ will be contained in $b_n\bar A^c$ for every $i=1,\ldots,N$, as soon as $n$ is large enough. Hence, for all large enough $n$,
$$\Big\{\sup_{1\leq i\leq N}\Big(\sup_{1\leq j\leq N} \langle X^{\e_n}_{t_j},\vf_i\rangle \Big) >0\Big\}
\subset \{X^{\e_n}\hbox{ hits }b_nA^c\}.$$
Using this inclusion and then (\ref{hittingtech2}) and (\ref{hittingtech1}), we immediately obtain (\ref{lowerboundhitting}).

We next turn to the more difficult upper bound
\begin{equation}
\label{upperboundhitting}
\limsup_{n\to\infty}\
\P_{[\e_n^{-1}]\delta_0}(X^{\e_n}{\rm\ hits\ }b_nA^c)\leq  P_{\delta_0}(Y^{(\kappa)}{\rm\ hits\ }A^c).
\end{equation}
We fix $\delta>0$ small enough so that the closed ball of radius $4\delta$ centered at $0$ is contained in $A$. As in Lemma \ref{analytic}, we let $A_\delta$ be the connected component of the open set
$$\{x\in A : {\rm dist}(x,A^c)>\delta\}$$
that contains $0$. We denote the exit measure from $A_\delta$ for the rescaled branching Brownian motion $X^{\e_n}$ by ${\cal E}^n_\delta$. In other words, the measure ${\mathcal{E}}^n_\delta$ is equal to $\e_n$ times the sum of the Dirac point masses at all points of $\partial A_\delta$ which are first exit points from $A_\delta$ for one of the historical paths associated with $X^{\e_n}$ (these historical paths are defined in Subsection \ref{subs: BBM et ob} for the branching Brownian motion $Z^{\e_n}$, and this definition is extended to $X^{\e_n}$ by an obvious
scaling transformation).

Let $\Phi$ be a continuous function on $\R^d$ such that $0\leq \Phi\leq 1$,
$\Phi=0$ on $A_{3\delta}$ and $\Phi=1$ on $A_{2\delta}^c$. Then, for
every $\eta>0$ and $\rho>0$,
\begin{eqnarray}
\label{upperhittingtech1}
\P_{[\e_n^{-1}]\delta_0}(X^{\e_n}{\rm\ hits\ }b_nA^c)
&=&\P_{[\e_n^{-1}]\delta_0}(X^{\e_n}{\rm\ hits\ }b_nA^c,\,
\langle {\cal E}^n_\delta,1 \rangle <\eta)\nonumber\\
&+&\P_{[\e_n^{-1}]\delta_0}\bigg(X^{\e_n}{\rm\ hits\ }b_nA^c,\,
\langle {\cal E}^n_\delta,1 \rangle \geq\eta, \,\int_0^\infty \langle X^{\e_n}_s,\Phi\rangle
ds\leq \rho\bigg)\nonumber\\
&+&\P_{[\e_n^{-1}]\delta_0}\bigg(X^{\e_n}{\rm\ hits\ }b_nA^c,\,
\langle {\cal E}^n_\delta,1 \rangle \geq\eta, \,\int_0^\infty \langle X^{\e_n}_s,\Phi\rangle
ds> \rho\bigg).
\end{eqnarray}
Let $\alpha_n(\eta)$, $\beta_n(\eta,\rho)$ and $\gamma_n(\eta,\rho)$ be the three terms appearing in the right-hand side of (\ref{upperhittingtech1}) in this order.

We first bound $\alpha_n(\eta)$. Provided $n$ is sufficiently large, $b_nA^c$ is contained in $A^c_{\delta/2}$ and thus
$$\alpha_n(\eta)\leq \P_{[\e_n^{-1}]\delta_0}(X^{\e_n}{\rm\ hits\ }A^c_{\delta/2},\, \langle {\cal E}^n_\delta,1 \rangle <\eta).$$
Note that the times at which the historical paths of $X^{\e_n}$ exit $A_\delta$ form a stopping line in the sense of \cite{Ch}. We can thus apply the strong Markov property at a stopping line (Proposition 2.1 in \cite{Ch}) to see that $\alpha_n(\delta)$ is bounded above by the probability for a branching Brownian motion (without killing) starting initially with less than $\eta \e_n^{-1}$ particles, that one of the historical paths reaches a distance greater than $\delta/(2\sqrt{\e_n})$ from its starting point (to be precise we need a slight extension of the results in \cite{Ch}, since our spatial motion is not standard Brownian motion, but Brownian motion killed inside $\Gamma_\varpi$). The estimate (\ref{cons-Sawyer}) now gives
\begin{equation}
\label{upperhittingtech2}
\alpha_n(\eta)\leq C''_1(d,\nu)\, \frac{4\eta}{\delta^2}.
\end{equation}

Then, we have 
$$\beta_n(\eta,\rho) \leq \P_{[\e_n^{-1}]\delta_0}\bigg(\langle {\cal E}^n_\delta,1 \rangle \geq\eta, \,\int_0^\infty \langle X^{\e_n}_s,\Phi\rangle ds\leq \rho\bigg).$$
Recall that $\Phi=1$ on $A^c_{2\delta}$ and in particular $\Phi=1$ on $\overline B(x,\delta)$ for every $x\in \partial A_\delta$. We use the strong Markov property at the same stopping line as in the previous argument, together with a simple coupling argument, to write that
$$\beta_n(\eta,\rho) \leq \P_{[\eta\e_n^{-1}]\delta_0}\bigg(\int_0^\infty \langle \widetilde X^{\e_n}_s,\ind_{\overline B(0,\delta)}\rangle ds\leq \rho\bigg),$$
where $\widetilde X^{\e_n}$ is defined in terms of a branching Brownian motion $\widetilde Z^{\e_n}$ in the same way as $X^{\e_n}$ was defined from $Z^{\e_n}$. This branching Brownian motion $\widetilde Z^{\e_n}$ has the same offspring distribution as $Z^{\e_n}$, but particles are now killed at rate $\e_n$ homogeneously over $\R^d$. Furthermore, $\widetilde Z^{\e_n}$ also starts from $k\delta_0$ under the probability measure $\P_{k\delta_0}$. 

\smallskip
By Proposition \ref{convBBM}, the law of $(\widetilde X^{\e_n}_t)_{t\geq 0}$  under $\P_{[\eta\e_n^{-1}]\delta_0}$ converges as $n\to\infty$ to the law of $Y^{(1)}$ under $P_{\eta\delta_0}$ (in the notation of Subsection \ref{subs: hitting prob}). Noting that, for every fixed $s>0$, $Y^{(1)}_s$ a.s.
does not charge the boundary of the ball $\overline B(0,\delta)$, it follows that
$$\limsup_{n\to\infty}  \P_{[\eta\e_n^{-1}]\delta_0}\bigg( \int_0^\infty \langle \widetilde X^{\e_n}_s,\ind_{\overline B(0,\delta)}\rangle ds\leq \rho\bigg)\leq P_{\eta\delta_0}\bigg(\int_0^\infty \langle Y^{(1)}_s,\ind_{\overline B(0,\delta)}\rangle ds\leq \rho\bigg)=: \beta_\infty(\eta,\rho).$$
The continuity of sample paths of $Y^{(1)}$ ensures that $\beta_\infty(\eta,\rho)\longrightarrow 0$ as $\rho\to 0$, for every fixed $\eta >0$.

For the term $\gamma_n(\eta,\rho)$, we simply use the bound
$$\gamma_n(\eta,\rho) \leq \P_{[\e_n^{-1}]\delta_0}\bigg(\int_0^\infty \langle X^{\e_n}_s,\Phi\rangle ds> \rho\bigg).$$
This bound and the weak convergence of Theorem \ref{theo cv to sbm} imply that 
$$\limsup_{n\to\infty} \gamma_n(\eta,\rho)\leq P_{\delta_0}\Big(\int_0^\infty 
\langle Y^{(\kappa)}_s,\Phi\rangle
ds\geq  \rho\Big) \leq P_{\delta_0}(Y^{(\kappa)}\hbox{ hits }A^c_{3\delta}),$$
since $\Phi=0$ on $A_{3\delta}$. (To justify the first inequality in the last display, we also use the
fact that the extinction times of $X^{\e_n}$ under $ \P_{[\e_n^{-1}]\delta_0}$ are stochastically bounded,
which follows from a standard result in the case without killing.)

To complete the argument, fix $\vartheta>0$. By Lemma \ref{analytic} (ii), we can choose $\delta>0$ sufficiently small so that
$$P_{\delta_0}(Y^{(\kappa)}\hbox{ hits }A^c_{3\delta})\leq P_{\delta_0}(Y^{(\kappa)}\hbox{ hits } A^c) +\frac{\vartheta}{3}.$$
From (\ref{upperhittingtech2}), we can then choose $\eta>0$ sufficiently small so that for all large $n$,
$$\alpha_n(\eta)\leq \frac{\vartheta}{3}.$$
Finally we choose $\rho>0$ such that $\beta_\infty(\eta,\rho)\leq \frac{\vartheta}{3}$. From (\ref{upperhittingtech1}) and the previous estimates, we obtain
$$\limsup_{n\to\infty}\ \P_{[\e_n^{-1}]\delta_0}(X^{\e_n}{\rm\ hits\ }b_nA^c) \leq P_{\delta_0}(Y^{(\kappa)}\hbox{ hits } A^c) + \vartheta,$$
and since $\vartheta$ was arbitrary this completes the proof of (\ref{upperboundhitting}) and Lemma \ref{convhittingpro}. $\hfill\square$

\smallskip
We still have to discuss the case $a=0$ in (\ref{keylimit}). So, let us consider two sequences $(\e_n)_{n\geq1}$ and $(R_n)_{n\geq 1}$ such that $\e_nR_n^2\to 0$. Let $a_0>0$ and $\e'_n=\e_n\vee (a_0R_n^{-2})$. Since $\e_n\leq \e'_n$, we have
$$\liminf_{n\to\infty} R_n^2\, \P_{\delta_0}(Z^{\e_n}\hbox{ hits }R_nA^c) \geq \liminf_{n\to\infty} R_n^2\, \P_{\delta_0}(Z^{\e'_n}\hbox{ hits }R_nA^c) =u_{(\kappa a_0)}(0),$$
by the case $a>0$. By Lemma \ref{analytic} (i), $u_{(\kappa a_0)}(0)$ can be made arbitrarily close to $u_{(0)}(0)$ when $a_0$ is small, and so 
$$\liminf_{n\to\infty} R_n^2\, \P_{\delta_0}(Z^{\e_n}\hbox{ hits }R_nA^c) \geq u_{(0)}(0).$$
To obtain the corresponding upper bound, a similar coupling argument shows that it suffices to consider the case when $\e_n=0$ for every $n$, that is when there is no killing inside the obstacles. Hence, consider the branching Brownian motion $Z^0=Z^{\varpi,0}$ (the notation is even more legitimate since $Z^{\varpi,0}$ does not depend on $\varpi$). For every $\rho>0$, define a rescaled version of $Z^0$ by setting
$$\langle \overline X^{(\rho)}_t,\vf\rangle  =\rho\int Z^0_{\rho^{-1}t}(dx)\,\vf(\rho^{1/2}x).$$
By Proposition \ref{convBBM}, the law of $(\overline X^{(\rho)}_t)_{t\geq 0}$ under $\P_{[\rho^{-1}]\delta_0}$ converges to the law of $Y^{(0)}$ under $P_{\delta_0}$ as $\rho$ tends to $0$. Set $\rho_n=R_n^{-2}$, in such a way that
\begin{equation}
\label{upperhitting4}
\{Z^0\hbox{ hits }R_nA^c\}=\{\overline X^{(\rho_n)}\hbox{ hits }A^c\}.
\end{equation}
A simplified version of the arguments of the proof of Lemma \ref{convhittingpro} shows that
$$\limsup_{n\to\infty} \P_{[\rho_n^{-1}]\delta_0}(\overline X^{(\rho_n)}\hbox{ hits }A^c) \leq P_{\delta_0}( Y^{(0)}\hbox{ hits }A^c)=1-\exp(-u_{(0)}(0)).$$
Arguing as in the first part of the proof of the theorem and using (\ref{upperhitting4}) yields
$$\limsup_{n\to\infty} R_n^2\, \P_{\delta_0}(Z^{\e_n}\hbox{ hits }R_nA^c)
\leq u_{(0)}(0),$$
which completes the proof of Theorem \ref{theo charging proba}.
$\hfill\square$

\medskip
\noindent{\bf Proof of Corollary \ref{escapepro}.} For every $r>0$,
$$\P_{n_\e\delta_0}(\sqrt{\e}\,R^{\varpi,\e}< r)= \Big(1 - \P_{\delta_0}(Z^{\varpi,\e}\hbox{ hits }B(0,\e^{-1/2}r)^c)\Big)^{n_\e}.$$
However, Theorem \ref{theo charging proba} shows that, $\bP(d\varpi)$ a.s.,
$$n_\e\,\P_{\delta_0}(Z^{\varpi,\e}\hbox{ hits }B(0,\e^{-1/2}r)^c) = \frac{\e n_\e}{r^2}\times \Big(\frac{r^2}{\e}\,
\P_{\delta_0}(Z^{\varpi,\e}\hbox{ hits }B(0,\e^{-1/2}r)^c)\Big)$$
converges to $\frac{b}{r^2} \,u^\circ_{\kappa r^2}(0)$ as $\e \to 0$. The desired 
result follows. $\hfill\square$

\bigskip
\noindent J.F. Le Gall, D\'epartement de math\'ematiques, Universit\'e Paris-Sud, 91405 ORSAY C\'edex, France

e-mail: {\tt jean-francois.legall@math.u-psud.fr}

\smallskip
\noindent A. V\'eber, DMA, Ecole normale sup\'erieure, 45 rue d'Ulm, 75230 PARIS C\'edex 05, France

e-mail: {\tt amandine.veber@ens.fr}

\end{document}